\documentclass[10pt,a4paper,twoside]{article}
\usepackage{a4,amsfonts}

\newcommand{\Z}{{\mathbf{Z}}}   
\newcommand{\Q}{{\mathbf{Q}}}
\def\qed{\hspace{\fill}\hbox{${\vcenter{\vbox{              
	\hrule height 0.4pt\hbox{\vrule width 0.4pt height 6pt
	\kern5pt\vrule width 0.4pt}\hrule height 0.4pt}}}$}}

\newtheorem{theorem}{Theorem}
\newtheorem{lemma}{Lemma}[section]

\newtheorem{definitionhelp}[lemma]{Definition}
\newtheorem{remarkhelp}[lemma]{Remark}
\newtheorem{questionhelp}[lemma]{Question}
\newtheorem{corollary}[lemma]{Corollary}
\newtheorem{proposition}[lemma]{Proposition}
\newenvironment{definition}{\begin{definitionhelp}\rm}{\end{definitionhelp}}
\newenvironment{remark}{\begin{remarkhelp}\rm}{\end{remarkhelp}}
\newenvironment{question}{\begin{questionhelp}\rm}{\end{questionhelp}}

\def\runninghead#1#2{\pagestyle{myheadings}
\markboth{{\protect\footnotesize\it{\quad #1}}\hfill}
{\hfill{\protect\footnotesize\it{#2\quad}}}}
\input epsf
\newcommand{\epsfig}[1]{\noindent
	\epsfbox{#1}}
\newcommand{\cepsfig}[1]{\par\smallskip\noindent
	\centerline{\epsfig{#1}}\smallskip}
\newcommand{\twocepsfig}[2]{\par\smallskip\noindent
	\centerline{\epsfig{#1}}\par\smallskip
	\centerline{\epsfig{#2}}\smallskip}
\newcommand{\threecepsfig}[3]{\par\smallskip\noindent
	\centerline{\epsfig{#1}}\par\smallskip
	\centerline{\epsfig{#2}}\par\smallskip
	\centerline{\epsfig{#3}}\smallskip}
\newcommand{\eqnfig}[2]{
	$$\begin{array}{c}
	\epsfig{#1}
	\end{array} 
	\eqno{\rm (#2)}
	$$}
\newcommand{\hcenterfig}[1]{
	$$
	\begin{array}{c}
	\epsfig{#1}
	\end{array} 
	$$}

\newcommand{\ti}{\hat}

\begin{document}
%
%
\runninghead{\quad Jan Kneissler}{Properties of Ladders \quad}
\title{{\large\it On Spaces of connected Graphs I}\\[0.2cm] Properties of Ladders}
\author{Jan A.~Kneissler}

\date{\quad} 

\maketitle
\thispagestyle{empty}
\vspace{6pt}
\begin{abstract}
{
\noindent
We examine spaces of connected tri-/univalent graphs subject to
local relations which are motivated by the theory of Vassiliev invariants.
It is shown that the behaviour of ladder-like subgraphs is strongly related to the 
parity of the number of rungs: there are similar relations for 
ladders of even and odd lengths, respectively.
Moreover, we prove that - under certain conditions - an even number of 
rungs may be transferred from one ladder to another.
}
\end{abstract}

\section{Introduction}
This note is the first in a series of three papers and its main goal
is to provide relations that will be employed by the following ones.
Nevertheless, we have the feeling that the results are interesting
enough to stand alone.
The objects we deal with are combinatorial multigraphs in which all
vertices have valency $1$ or $3$, equipped with some additional data: 
at each trivalent vertex a cyclic ordering of the three incoming edges
is specified and every univalent vertex carries a colour.
In pictures we assume a counter-clockwise ordering
at every trivalent vertex and indicate the colour of univalent vertices by integers.
The graphs have to be connected and must contain at least one trivalent
vertex, but we do not insist that there are univalent vertices.
For simplicity, such graphs will be called {\sl diagrams} from now on.

The motivation is given by the combinatorial access to Vassiliev invariants.
We briefly mention the most relevant facts about this comprehensive 
family of link invariants. Vassiliev invariants with values in
a ring $R$ form a filtered algebra. A deep theorem of Kontsevich
(see \cite{Ko} and \cite{BN}) states that for $\Q\subset R$ the filtration quotients are isomorphic
to the graded dual of the Hopf-algebra of chord diagrams. By the structure
theorem of Hopf-algebras it suffices to examine the primitive elements.

The algebra of chord diagrams is rationally isomorphic (\cite{BN})
to an algebra of tri-/univalent graphs (sometimes called
Chinese characters), defined by relations named (AS) and (IHX) (see 
Definition \ref{defB}). 
The coproduct is given by distributing the components into two groups,
so the primitives are spanned by connected graphs.
The colours of the univalent vertices represent the link-components 
and half the total number of vertices corresponds to the degree
of Vassiliev invariants. The number of univalent vertices allows
a second grading that corresponds to the eigenspace-decomposition 
of primitive Vassiliev invariants with respect to the cabling operation.

The most important open question in Vassiliev theory is whether all
invariants of finite degree taken together form a complete invariant
of knots and links. For knots there is a weaker (but equally essential)
question: 

\begin{question}
\label{questinv}
Are there any (rational valued) Vassiliev invariants that are able to detect non-invertibility
of knots?
\end{question}

\noindent
There is a weak hope that this question might be settled in the combinatorial
setting, where it translates into the question whether all 
diagrams with an odd number of univalent vertices vanish.

A progress in the combinatorics of connected graphs (\cite{Vo}) was initiated by
the observation that replacing a trivalent vertex by the graph 
that is shown in relation (x) in the next section is a well-defined operation
(i.e.~independent of the choice of a trivalent vertex and the orientation).
This operation called $x_n$ is an element of an algebra called $\Lambda$ 
that acts on spaces of connected graphs. 
One purpose of our investigations (and in fact the initial one) was
to obtain a family of relations in $\Lambda$. The result is presented in \cite{CC2}.

Our philosophy is to declare diagrams that are coming from lower degrees 
by the action of $\Lambda$ as uninteresting, which explains why we factor
by the relation (x) in $\ti B^u$ (see Definition \ref{defB}).
This can be justified in two ways: First, note that Question
\ref{questinv} is equivalent to asking whether $\ti B^u = 0$ for all odd
$u$.
Second, any good upper bound for the dimensions of $\ti B^u$ will lead
to good upper bounds for $\dim B^u$ as we demonstrate in \cite{CC3}.     
 
\section{Results}
First let us introduce notations for certain subgraphs of diagrams. 
An edge that connects a trivalent and a univalent vertex is called {\sl leg}.
For $n \geq 2$, a subgraph consisting of $3n+2$ edges of the following type 
\cepsfig{}
\noindent
is called a $n$-{\sl ladder}. The two uppermost and the two lowermost edges are called 
{\sl ends} of the ladder. They may be connected to univalent or trivalent vertices
in the rest of the diagram. The $3n-2$ other edges are called 
{\sl interior} of the ladder. 
A ladder is said to be {\sl odd} or {\sl even} according to the parity of $n$.
A $2$-ladder is also called {\sl square}. 
Finally, for a partition $u = (u_1,\ldots,u_n)$, $u_1 \leq u_2 \leq \cdots u_n$,
we say a diagram is $u$-coloured if its univalent vertices carry the colours
$1$ to $n$ and there are exactly $u_i$ univalent vertices of colour $i$
for $1 \leq i \leq n$.

\begin{definition}
\label{defB}
For any partition $u$ let
$$B(u) \;\;:=\;\; \Q\langle \;u\mbox{-coloured diagrams }\rangle\mbox{ } / 
\mbox{ (AS), (IHX)}\;\;\mbox{ and } \;\; \ti B(u) \;\;:=\;\; B(u)\mbox{ } / \mbox{ (x),}$$
where (AS), (IHX) and (x) are the following local relations:
\eqnfig{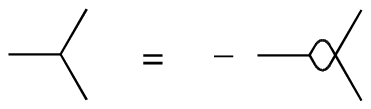}{AS}
\eqnfig{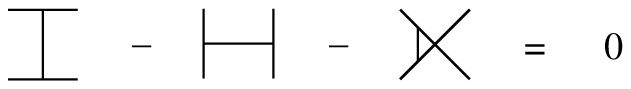}{IHX}
\eqnfig{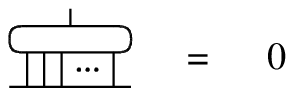}{x}
The extremal partitions $u_1 = \cdots = u_n = 1$ and $n=1$ are most
important, so they get their own names:
$F(n) := B\big((1,\cdots,1)\big)$, $\ti F(n) := \ti B\big((1,\cdots,1)\big)$ and
$B^u := B\big((u)\big)$, $\ti B^u := \ti B\big((u)\big)$.
\end{definition}

\begin{remark}
The careful reader will notice that in our calculations we have to divide by $2$ and
$3$ only, so all statements remain valid if one works with $\Z[\frac16]$-modules
instead of $\Q$-vectorspaces. 
\end{remark}

\noindent
Let us mention an important consequence of the relations (AS), (IHX) and (x):
\eqnfig{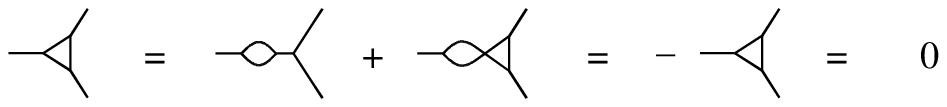}{t}

\noindent
By (IHX) and (t) the ends of a square may be permuted:
\cepsfig{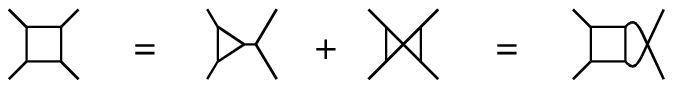}

\noindent 
This obviously implies that neighbouring ends of any $n$-ladder can be interchanged.
But it is a little surprise that for all even ladders any permutation of the ends 
yields the same element, i.e.~for $(2m)$-ladders ($m \geq 1$) we have the following relation called LS (``ladder
symmetry''): 
\eqnfig{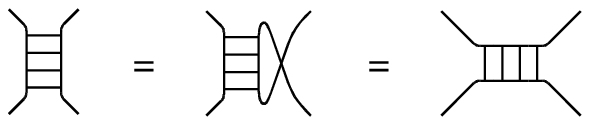}{LS}

\noindent
Here and later on, when making statement about $n$-ladders, we draw the corresponding 
picture for a generic but small value of $n$. This should not lead to confusion.

\noindent
For any $m \geq 1$, we have the following relation of the
IHX-type, involving odd ladders of length $2m+1$:
\eqnfig{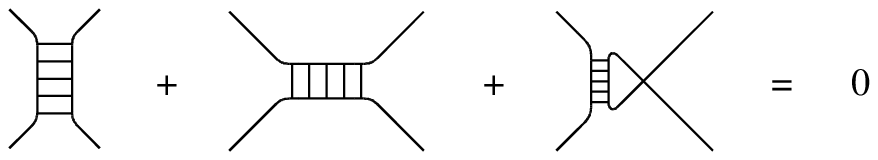}{LIHX}

\noindent
We present a further relation named LI, in which an edge has been 
glued to non-neighbouring ends of an odd ladder:
\eqnfig{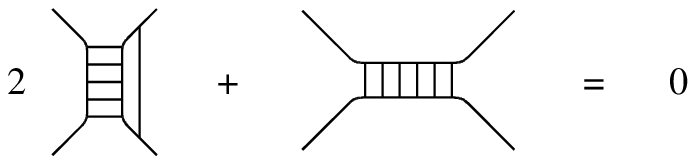}{LI}

\noindent
Finally, we have the following relation LL that replaces two parallel
ladders of lengths $2n+1$ and $2m+1$ 
by a single $(2n+2m+2)$-ladder:
\eqnfig{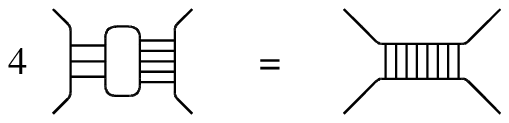}{LL}

\begin{theorem}
\label{theo1}
The relations (LS), (LIHX), (LI) and (LL) are valid in $\ti B(u)$ for any partition $u$.
\end{theorem}

\begin{definition}
{\rm
Suppose we have a diagram $D$ containing two ladders $L_1$, $L_2$ 
which have disjoint interiors. If we remove the interiors 
of $L_1$ and $L_2$ from $D$ we get a $-$ possibly disconnected $-$
graph $D^\prime$. Let $D^\prime_1, \ldots D^\prime_k$ denote 
the components of $D^\prime$ that contain at least one end of $L_1$ {\it and}
at least one end of $L_2$.
If $k=1$ and $D^\prime_1$ is a tree and 
$D^\prime_1$ contains {\it exactly} one end of $L_1$ and {\it exactly} one end of $L_2$, then we say that $L_1$ and $L_2$ 
are {\sl weakly connected}. Otherwise we call $L_1$ and $L_2$ 
{\sl strongly connected}.
If the intersection of the interiors of $L_1$ and $L_2$ is not empty
(which implies that $L_1$ and $L_2$ are sub-ladders of a single longer
ladder) then $L_1$ and $L_2$ are also called strongly connected.
}
\end{definition}

\noindent
Example: In the following diagrams 
the $3$-ladders and $4$-ladders are strongly connected to the square but 
weakly connected to each other (so being strongly connected is not a 
transitive relation):
\cepsfig{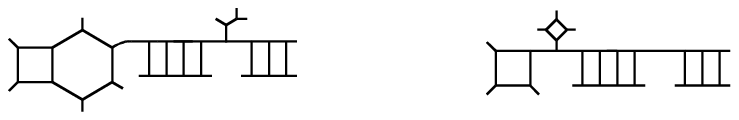}

\noindent
For the rest of this section let us assume that $D$ is a $u$-coloured 
diagram with two specified $3$-ladders $L_1$ and $L_2$. 
For $a,b \geq 2$ let $D_{a,b}$  denote the diagram that is obtained by replacing $L_1$ and $L_2$ by two
ladders (in the same orientation) of length $a$ and $b$, respectively.

\begin{theorem}[Square-Tunnelling relation]
\label{theo2}
If $L_1$ and $L_2$ are strongly connected in $D$ then $D_{2,4} \;=\; D_{4,2}$ 
in $\ti B(u)$. 
\end{theorem}

\noindent
If either $L_1$ or $L_2$ is subladder of a longer ladder, then $L_1$ or $L_2$ are automatically
strongly connected, which allows to make a nice statement:

\begin{corollary}
For any $a \geq 2$, $b \geq 4$ with $a+b \geq 7$ 
we have $D_{a,b} \;=\; D_{a+2,b-2}$ in $\ti B(u)$.
\end{corollary}

\noindent
One might ask whether the strong connectivity condition is essential
in Theorem \ref{theo2}, 
in other words:

\begin{question}
Are there any $u$-coloured diagrams $D$ with {\it weakly} connected ladders 
such that $D_{2,4} \neq D_{4,2}$? 
\end{question}

\noindent
At least for the case length$(u)=1$, we may give the answer:
\begin{theorem}[Square-Tunnelling relation in $\ti B^ u$]
\label{coro}
Let $D$ be a diagram of $B^u$ (i.e.~all univalent vertices carry the same colour),
then $D_{2,4} = D_{4,2}$. 
\end{theorem}

\noindent
Compared to the other relations, which are local, the square-tunnelling 
relation has a completely different character: it relates ladders that
might be located arbitrarily far apart in a diagram.
We already mentioned in the introduction a situation
where a global structure (the action of $\Lambda$ on $B(u)$) emerges from 
local relations ((IHX) and (AS)).
There, it is easy to understand how the subgraphs $x_n$ (shown in relation (x))
move around, because they can go from one vertex to a neighbouring one. 
But here, we are not able to trace the way of the square from the ladder 
in which it disappears to the other ladder.
It is just like the quantum-mechanical 
effect where an electron tunnels through a classically impenetrable 
potential barrier: one can calculate that it is able to go from one place 
to the other, but we cannot tell how it actually {\it does} it.

\section{The local ladder relations}

Let us start with two well-known and frequently-used implications of
(IHX), (AS) and (t).
We assume that we are given a diagram $D$ of $\ti F(n)$ together with an arbitrary
grouping of the legs of $D$ into two classes called ``entries'' and ``exits''.
Say there are $p$ entries and $q$ exits (thus $p+q=n$).
For $1 \leq i \leq p$ ($1\leq j \leq q$) let $D_i$ ($D^j$)  denote the elements
of $\ti F(n+1)$ one obtains when the $(n+1)$-th leg is glued to the $i$-th entry 
(the $j$-th exit) of $D$. 
Furthermore for $1 \leq i < j \leq p$ ($1\leq k < l \leq q$) 
let $D_{ij}$ ($D^{kl}$) denote the elements
of $\ti F(n)$ having an additional edge between
the $i$-th and $j$-th entry (the $k$-th and $l$-th exit) of $D$.
$D_i$, $D^j$, $D_{ij}$, $D^{kl}$ look typically like this:
\cepsfig{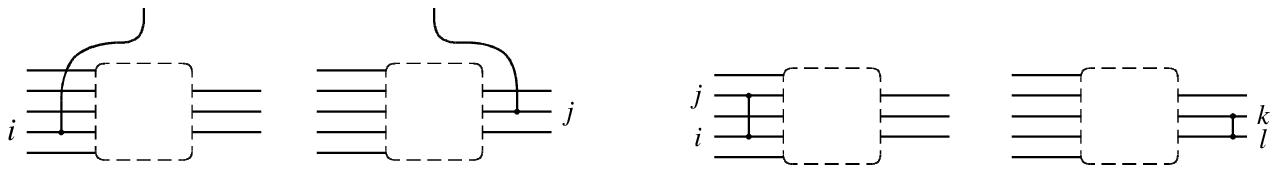}

\begin{lemma}
\label{lemmab}
For any $D \in \ti F(n)$ with $p$ entries and $q$ exits, we have 
\begin{itemize}
\item[a)]
\quad$\displaystyle\sum_{i=1}^{p} D_i \;\;=\;\; \sum_{j=1}^{q} D^j$
\item[b)]
\quad$\displaystyle\sum_{1 \leq i < j \leq p} D_{ij} \;\;=\;\; \sum_{1 \leq k < l \leq q} D^{kl}$.
\end{itemize}
\end{lemma}
\begin{proof}
We cut the box into little slices, such that each slice is of one of
the following types:
\cepsfig{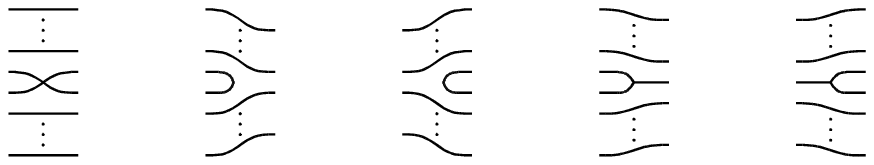}

\noindent
It is easy to verify that a) and b) are valid for slices of these types. 
\end{proof}
\qed
\bigskip

\newcommand{\x}{c_2}
\newcommand{\y}{c_3}
\newcommand{\z}{c_1}

\noindent
Let us introduce a notation for certain elements of 
$\ti F(4)$. For a word $w$ in the letters $\z$, $\x$, $\y$, let $\langle w\rangle$
be a diagram of $\ti F(4)$ that is constructed as follows: Take three strings,
put edges between the strings according to $w$, and then glue 
string $1$ and string $2$ together on the right side. $\z$, $\x$ and $\y$ correspond to
$1$-$2$, $2$-$3$ and $1$-$3$ edges, respectively. For example 
$\langle \x\y\z^2\y\rangle$, $\langle \x^4\rangle$, $\langle \y^4\rangle$
are the following diagrams:
\cepsfig{}

\begin{proposition}
\label{propxy}
For $m \geq 1$ and any word $u$ in the letters $\z$, $\x$, $\y$, we 
have $\langle u\x^{2m}\rangle \;\;=\;\; \langle u\y^{2m}\rangle$.
\end{proposition}
\begin{proof}
The statement is implicitly hidden in section 5 of \cite{Vo}, but to the readers 
convenience we will give a simple inductive proof of the statement. 
We know already that it is true for $m = 1$. 
We have the following equalities for arbitrary words $u$, $v$:

$$\langle u\z v \rangle + \langle u\x v \rangle + \langle u\y v \rangle \;\;=\;\; 0
\eqno{(\mbox{i})}$$
$$ \langle u\z^{}\x^n \rangle   \;\;=\;\;  \langle u\z^{}\y^n \rangle\eqno{(\mbox{ii})}$$
$$ \langle u\x^{}\z^{}\x^n \rangle  \;\;=\;\;  \langle u\x^{}\y^{}\x^n \rangle \mbox{\quad for }n\geq 2\eqno{(\mbox{iii})}$$
$$ \langle u\y^{}\z^{}\y^n \rangle  \;\;=\;\;  \langle u\y^{}\x^{}\y^n \rangle \mbox{\quad for }n\geq 2\eqno{(\mbox{iv})}$$
(i) is due to Lemma \ref{lemmab} b). (ii) and (iii) are implications of
(IHX) and (x):
\twocepsfig{}{}

\noindent
(iv) is shown similarly to (iii), interchanging string $1$ with string $2$.
(i), (iii) and (iv) imply 
$$2\langle u\x^{}\y^{}\x^n\rangle \;\;=\;\; -\langle u\x^{n+2}\rangle
\mbox{\quad and \quad}
2\langle u\y^{}\x^{}\y^n\rangle \;\;=\;\; -\langle u\y^{n+2}\rangle
\mbox{\quad for } n \geq 2.
\eqno{(\mbox{v})}
$$
Assuming $m \geq 2$ and $\langle w\x^{2m-2}\rangle \;\;=\;\; \langle w\y^{2m-2}\rangle$ by induction hypothesis, 
we obtain
\begin{eqnarray*}
\langle u\x^{2m} \rangle  - \langle u\y^{2m} \rangle  
&\stackrel{\mbox{\tiny(i)}}{=}& 
-\langle u\y^{}\x^{2m-1} \rangle -\langle u\z^{}\x^{2m-1} \rangle +\langle u\x^{}\y^{2m-1} \rangle  +\langle u\z^{}\y^{2m-1} \rangle  \\
&\stackrel{\mbox{\tiny(ii)}}{=}& -\langle u\y^{}\x^{2m-1} \rangle +\langle u\x^{}\y^{2m-1} \rangle \\ 
&\stackrel{\mbox{\tiny I.H.}}{=}&
-\langle u\y^{}\x^{}\y^{2m-2} \rangle +\langle u\x^{}\y^{}\x^{2m-2} \rangle  
\;\;\stackrel{\mbox{\tiny(v)}}{=}\;\; \frac12\langle u\y^{2m} \rangle -\frac12\langle u\x^{2m} \rangle.
\end{eqnarray*}
\qed
\end{proof}
\bigskip

\noindent
{\bf Proof of Theorem \ref{theo1}}\quad
The first equality of the (LS)-relation is true by Proposition \ref{propxy}. 
To prove the second equality, we simply swap the upper ends
of the ladder and use the first equality of (LS):
\cepsfig{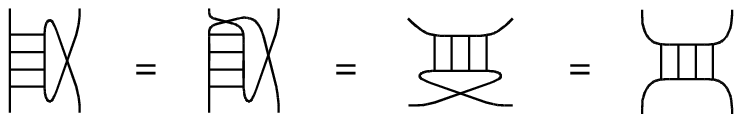}

\noindent
If we take three ends of a $(2m)$-ladder as entries and the forth as exit, apply 
Lemma \ref{lemmab} b) and (LS), we get the (LIHX)-relation: 
\cepsfig{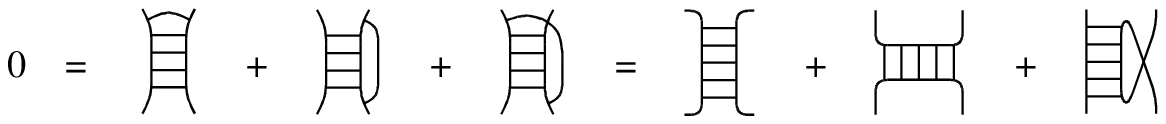}

\noindent
Relation (LI) is obtained when we append an edge between the 
ends on the right side of (LIHX) and use (IHX) and (x):
\cepsfig{}

\noindent
(LL) is the hardest one; we apply (IHX) and use Lemma \ref{lemmab} a):
\cepsfig{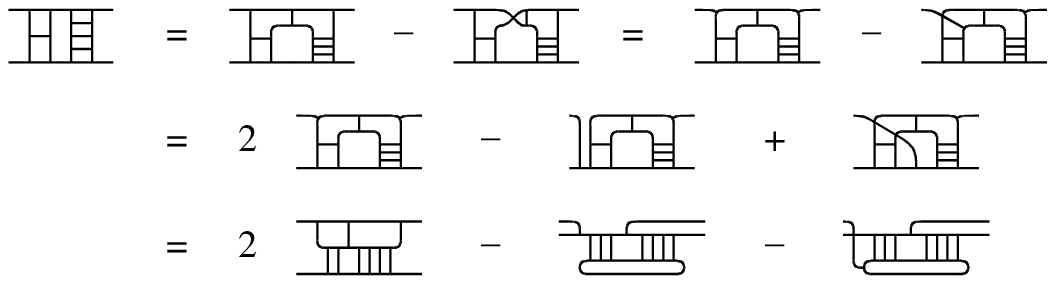}

\noindent
To each of the three resulting diagrams we apply the following identity
\cepsfig{}

\noindent
to obtain by (t), (IHX), (x), (LS) and (LI) the promised result:
\cepsfig{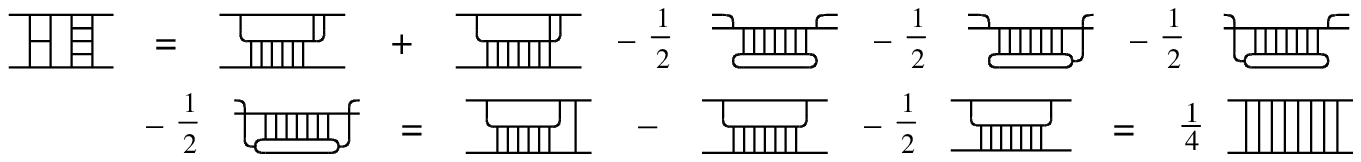}

\qed

\section{Relations in $\ti F(6)$}

\newcommand{\e}{{\;\;=\;\;}} 
\newcommand{\eq}[1]{{\;\;\stackrel{\mbox{\tiny(\ref{#1})}}{=}\;\;}} 
\newcommand{\eqq}[1]{{\stackrel{\mbox{\tiny(\ref{#1})}}{=}}} 
\newcommand{\eqz}[2]{{\;\;\stackrel{\mbox{\tiny(\ref{#1},\ref{#2})}}{=}\;\;}} 
\newcommand{\si}[1]{{\;\;\stackrel{\mbox{\tiny(\ref{#1})}}{\sim}\;\;}} 
\newcommand{\sinoref}[1]{{\;\;\stackrel{\mbox{\tiny({#1})}}{\sim}\;\;}}

\renewcommand{\c}[1]{c_{#1}^{}}
\newcommand{\C}[2]{c_{#1}^{#2}}
\renewcommand{\d}[1]{d_{#1}^{}}
\newcommand{\D}[2]{d_{#1}^{#2}}
\newcommand{\s}[1]{s_{#1}^{}}
\renewcommand{\S}[2]{s_{#1}^{#2}}
\renewcommand{\l}[1]{l_{#1}^{}}
\renewcommand{\_}[1]{_{#1}^{}}
\renewcommand{\H}[1]{H_{#1}^{}}
\newcommand{\sym}[3]{[#1]_{#2}^{#3}}
\newcommand{\symo}[1]{[#1]}

The key to Theorem \ref{theo1} is a certain equation in $\ti F(6)$. 
Written in the notation we will introduce in this section, it appears
misleadingly simple:
$$ \sym{y}13 \;\;=\;\; 0.$$
Nevertheless, it requires a lot of calculation. We have spent a considerable
amount of time trying to take the computation into a bearable form. 
We hope $-$ for the readers sake $-$ that this has not been a vain endeavour.

\subsection{The algebra $\Xi$}
\label{sectxi}

We do not like to draw pictures all the time, so we introduce an algebra 
that enables us to write down sufficiently many elements of $\ti F(6)$
in terms of a few number of symbols.

Let $\Xi$ be the space of tri-/univalent graphs with exactly six numbered
univalent vertices quotiented by the (AS), (IHX) and (x). 
The only difference between $\Xi$ and $\ti F(6)$ is that in $\Xi$ 
we do not require the graphs to be connected. Also we do not insist
that there actually {\it are} any trivalent vertices.

All graphs in an (AS)- or (IHX)-relation have same number of components. 
This implies that $\Xi = \ti F(6) \oplus \Delta$ where $\Delta$ is the subspace of
$\Xi$ that is spanned by disconnected graphs. So we can project any relation
that we find in $\Xi$ to $\ti F(6)$. 

Let us introduce a multiplication in $\Xi$: For any two graphs $D_1, D_2$
let $D_1D_2$ denote the graph that is obtained when the univalent vertices 
$6$, $5$, $4$ of $D_1$ are glued to univalent vertices  $1$, $2$, $3$ of
$D_2$, respectively. This extends to a well defined map $\Xi \otimes \Xi \rightarrow \Xi$, so
we may view $\Xi$ as algebra. It is a graded algebra, but we make no explicit use of this.
In pictures we always arrange the ends always like this:
\cepsfig{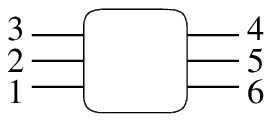}

\noindent
Let $\c1,\,\c2,\,\c3,\,\s1,\,\s2,\,\d1,\,\d2,\,y$ denote the elements of $\Xi$ that are represented by the following graphs:
\cepsfig{}

\noindent
Furthermore let $\l1$, $\l2$, $\l3$, $\l4, \ldots$ be the family of
elements of $\Xi$ that is represented by the following graphs
($l_n$ contains a $n$-ladder for $n \geq 2$ and $\l1 = \c1\c2$):
\cepsfig{}

\noindent
Let $z := \c1+\c2+\c3$, then by Lemma \ref{lemmab} a), $z$ is central in $\Xi$. 
We will continually substitute $\c3$ by $z-\c1-\c2$ and collect all $z$'s 
at the beginnings of words.

\begin{proposition} 
\label{proprel}
The following relations are fulfilled in $\Xi$:
\begin{eqnarray}
\nonumber
s_1^2 &=& s_2^2 \e (\s1\s2)^3 \e 1 \\
\nonumber
d\_i c\_i &=& {\smash{c\_id\_i \e 0}} \mbox{\quad for }i = 1,2 \\
\nonumber
s\_ic\_i &=& \smash{c\_is\_i} \mbox{\quad for }i = 1,2 \\
\label{screl}
s\_ic\_j &=& c\_ks\_i \mbox{\quad for }\{i,\,j,\,k\} = \{1,\,2,\,3\} \mbox{ and }i \neq 3\;\;\;\; \\
\label{cdsrel}
d_i^{} &=& c\_i\,(1-s\_i) \mbox{\quad for }i = 1,2 \\
\label{cdcrel}
\c1\d2\c1 &=& \c2\d1\c2  \\
\label{dcccdrel}
\d2c_1^{2n+1}\d2 &=& \d2\c1c_2^{2n-1}\c1\d2 \e z\d2c_1^{2n}\d2
\e -\frac12z^{2n+2}\d2 \mbox{\quad for }n \geq 1 \\
\label{lsrel}
l\_n\,(1+\s1) &=& \d1\c2c_1^{n-1} + \c3c_1^n \mbox{\quad for }n \geq 2 \\
\label{slrel}
\s2\,l\_n &=& -l\_n + c_2^{n-1}\c1\d2 + c_2^n\c3 \mbox{\quad for }n \geq 2 \\
\label{yrel}
\c1y\c1 &=& \s1\c1y\c1 \e z\c1\d2\c1 - z\c1\s1\s2\d1\c2 \\
\label{zzssdcrel}
z^2\d2\c1 &=& z^2\s1\s2\d1\c2
\end{eqnarray}
\end{proposition}

\noindent
\begin{proof}
The first six equalities are obvious: $\d{i}\c{i}$ and $\c{i}\d{i}$ contain triangles, 
$\s{i}\c{i} = \c{i}\s{i}$ is a double application of (AS)
and $\d{i} = \c{i} - \s{i}\c{i}$ is just (IHX).
$\c1\d2\c1$ and $\c2\d1\c2$ are identical graphs:
\cepsfig{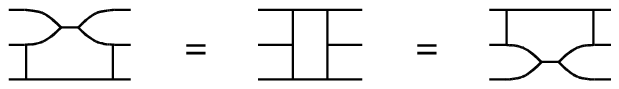}

\noindent
We will use in the sequel that due to Lemma \ref{lemmab} b), $z^k\d2$ is represented by
the following graph:
\cepsfig{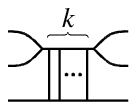}

\noindent
Then all equalities in (\ref{dcccdrel}) are simple consequences of (LI) and (LS):
\threecepsfig{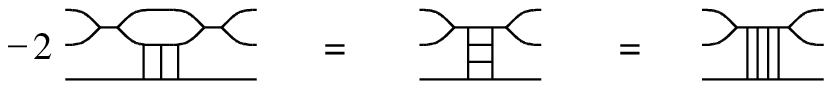}{}{}

\noindent
To show (\ref{lsrel}) and (\ref{slrel}), we have to use Lemma \ref{lemmab} a)
(for (\ref{slrel}) simply rotate the picture by $\pi$):
\cepsfig{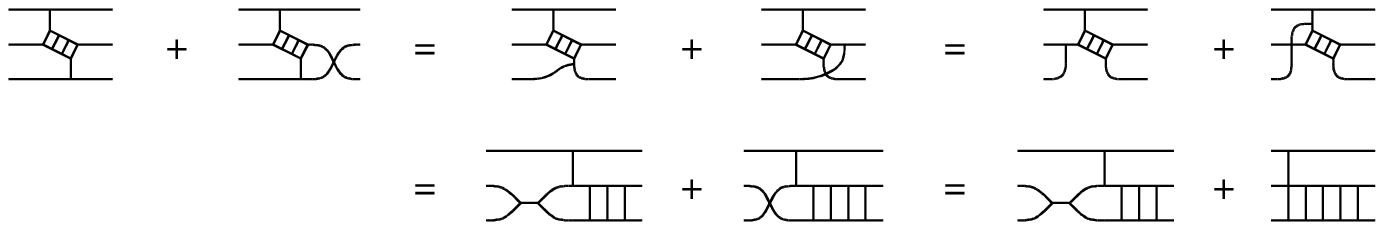}

\noindent
Finally we use (t) and (IHX) to show (\ref{yrel}) and (\ref{zzssdcrel}):
\twocepsfig{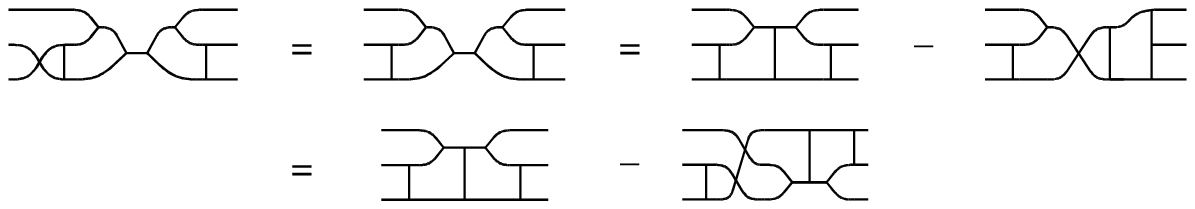}{}
\qed
\end{proof}

\noindent
Let us elaborate some consequences of the relations, which we will use later.
(Note that the application of a previous equation is indicated by stacking 
its number atop the equality sign but we do not mention the use
of equations without number.)

\begin{eqnarray}
\label{sdrel}
s\_id\_i & = & d\_is\_i \eq{cdsrel} s\_i(1-s\_i)c\_i \e (s\_i-1)c\_i \eq{cdsrel} -d\_i \\
\label{ddrel}
d_i^2 & \eqq{cdsrel} & d\_ic\_i - d\_is\_ic\_i \e 2c\_id\_i \e 0 \\
\label{dcdrel}
\d2\c1\d2 & \eqq{sdrel}& \frac12\d2\c1\d2 + \frac12\d2\s2\c1\s2\d2 \eq{screl} \frac12\d2(\c1+\c3)\d2  
\e  \frac12\d2(z-\c2)\d2 \eq{ddrel} 0 \\
\label{ddcrel}
\d2\d1\c2 & \eqq{cdsrel} & (1-\s2)\c2\d1\c2 \eq{cdcrel} (1-\s2)\c1\d2\c1 \eqz{screl}{sdrel} 
\c1\d2\c1 + \c3\d2\c1 \e z\d2\c1 \\
\label{dddrel}
\d2\d1\d2 & \eqq{cdsrel} & \d2\d1\c2(1-\s2) \eq{ddcrel} z\d2\c1(1-\s2) 
\eqz{screl}{sdrel} z\d2(\c1+\c3) \e  z^2\d2 
\end{eqnarray}
\begin{eqnarray}
c_2^3 + 2\c2\c1\c2 -zc_2^2 &=& \c2(2\c1+\c2-z)\c2 \e \c2\c1\c2 - \c2\c3\c2  \eq{screl} 
\c2\c1\c2 - \c2\s2\c1\s2\c2  \nonumber \\  
\label{cccrel}
& \eqq{cdsrel}  & \c2\c1\c2-(\c2-\d2)\c1(\c2-\d2) \eq{dcdrel} \d2\c1\c2 + \c2\c1\d2 
\end{eqnarray}
\begin{eqnarray}
\nonumber
z^2\d2(z-\c1) & = & z^2\d2\c3 \eq{dddrel} \d2\d1\d2\c3 \eqz{cdsrel}{dcdrel} -\d2\s1\c1\d2\c3
\e -\d2\s1(z-\c2-\c3)\d2\c3 \\
\nonumber
&\eqq{screl}& -z\d2\s1\d2\c3 + \d2\c2\s1\d2\c3  \eq{sdrel} z\d2\s2\s1\d2\c3
\eqz{cdsrel}{screl} z\d2\d1\s2\s1\c3 \\
\label{zzdcrel}
&\eq{screl}& z\d2\d1\s2\c2\s1 \eq{cdsrel}
z\d2\d1(\c2-\d2)\s1 \eqz{ddcrel}{dddrel} z^2\d2\c1\s1 - z^3\d2\s1  \\[0.4cm]
\label{cccdrel}
\nonumber
(\c2-\d2)c_1^2\d2 & \eq{cdsrel} & \c2\s2c_1^2\d2 \eqz{screl}{sdrel} 
- \c2c_3^2\d2 \e - \c2(z-\c1-\c2)^2\d2 \\
& = & 2z\c2\c1\d2 -c_2^2\c1\d2 -\c2c_1^2\d2 
\end{eqnarray}

\subsection{The Subspace $[\Xi_0]$}
The symmetric group $S_6$ operates on $\Xi$ by permutation of the univalent 
vertices. This allows us to regard $\Xi$ as a $\Z[S_6]$-module. To prevent
confusion, we use a dot to indicate this operation (so $\sigma \in \Z[S_6]$ 
applied to $\xi\in \Xi$ will be written as $\sigma\cdot\xi$).
Note that in general $\sigma\cdot\xi_1\xi_2 \neq (\sigma\cdot\xi_1)\xi_2$.

The elementary transpositions in $S_6$ will be named $\tau_i := (i\;\;i\!+\!1)$.
Obviously, $\tau_1\cdot\xi = s_1\xi$, $\tau_2\cdot\xi = s_2\xi$, 
$\tau_4\cdot\xi = \xi s_2$ and $\tau_5\cdot\xi = \xi s_1$.
Furthermore let $\mu := (1\;\;6)(2\;\;5)(3\;\;4)$, then 
$\mu$ operates by mirroring along the $y$-axis (there is always an even number of
trivalent vertices, so (AS) causes no change of sign).
For abbreviation we  
introduce the symbol $\sym{*}{*}{*}$:  
$$\sym{\xi}ij \;\;:=\;\; (1-\mu)(1+\tau_1)(1+\tau_3)(1+\tau_5)\cdot(c_1^i\xi c_1^j) \mbox{\quad for any }
\xi \in \Xi \mbox{ and } i,j \geq 0.$$
The maps $\xi \rightarrow \sym{\xi}ij$ are well-defined $\Q$-endomorphisms of $\Xi$.
Obviously $\sym{\c1\xi}ij = \sym{\xi}{i+1}j$ and $\sym{\xi\c1}ij = \sym{\xi}i{j+1}$.
A useful feature of this notation is that whenever $\mu\cdot w = w$ then 
$\sym{z^kw}ii = 0$. For instance, this is the case if $w$ is a palindrome 
in the letters $s_*^{},\,c_*^{},\,d_*^{}$. Other pleasing properties
of $\sym***$ are given in the following proposition.

\noindent
\begin{proposition}
\label{propsym}
For any $\xi\in\Xi$, $i,j \geq 0$ we have
\begin{eqnarray}
\label{symrel1}
\sym{\xi}ij \e -\sym{\mu\cdot\xi}ji \;\;=  &\sym{\s1\xi}ij&  =\;\; \sym{\tau_3\cdot\xi}ij \e \sym{\xi\s1}ij\\
\label{symrel2}
\c1\sym{\xi}ij\c1 &=& \sym{\xi}{i+1}{j+1}\\
\label{symrel3}
\c1\sym{\xi}ij + \sym{\xi}ij\c1  &=& \sym{\xi}{i+1}{j} + \sym{\xi}i{j+1}.
\end{eqnarray}
Moreover, if $\mu\cdot\xi = \xi$ and $\sym{\xi}01 = 0$, then $\sym{\xi}ij = 0$
for all $i,\,j \geq 0$.
\end{proposition}
\begin{proof}
The last three factors of $\pi := (1-\mu)(1+\tau_1)(1+\tau_3)(1+\tau_5)$ commute 
with each other and their product commutes with $(1-\mu)$. Since
$\tau_i^2 = \mu^2 = 1$, we have $\pi = -\pi\mu = \pi\tau_1=\pi\tau_3=\pi\tau_5$, which 
implies (\ref{symrel1}). Equality (\ref{symrel2}) becomes clear if one realises that for $i \in \{1,\,3,\,5\}$    
$$
\begin{array}{rclcrcl}
\c1(\tau_i\cdot\xi) &=& \tau_i\cdot\c1\xi, && (\tau_i\cdot\xi)\c1 &=& \tau_i\cdot\xi\c1\\
\c1(\mu\cdot\xi) &=& \mu\cdot\xi\c1, && (\mu\cdot\xi)\c1 &=& \mu\cdot\c1\xi.
\end{array}$$
Let $x := (1+\tau_1)(1+\tau_3)(1+\tau_5)\cdot c_1^i\xi c_1^j$, then $\sym{\xi}ij = x-\mu\cdot x$ and
\begin{eqnarray*}
\c1\sym{\xi}ij + \sym{\xi}ij\c1 &=& \c1x - \c1(\mu\cdot x) + x\c1 - (\mu\cdot x)\c1
\e \c1 x - \mu\cdot x\c1 + x\c1 -\mu\cdot \c1x \\
&=& (1-\mu)\cdot\c1 x + (1-\mu)\cdot x\c1 \e \sym{\xi}{i+1}{j} + \sym{\xi}i{j+1}. 
\end{eqnarray*}
Finally, the following identity shows by induction that
$\sym{\xi}00 = \sym{\xi}01 = 0$ implies $\sym{\xi}0j = 0$ for $j \geq 2$:
$$
\sym{\xi}0{j} \e \sym{\xi}0{j} + \sym{\xi}1{j-1} - \sym{\xi}1{j-1}
\eqz{symrel3}{symrel2} \c1\sym{\xi}0{j-1} + \sym{\xi}0{j-1}\c1 - \c1\sym{\xi}0{j-2}\c1
$$
By (\ref{symrel2}) we obtain $\sym{\xi}ij = c_1^i\sym{\xi}0{j-i}c_1^i = 0$
for $i \leq j$. If $i > j$ we have $\sym{\xi}ij \eq{symrel1} -\sym{\mu\cdot\xi}ji = \sym{\xi}ji = 0$.  
\end{proof}
\qed

\noindent
\begin{definition}
Let $\Xi_0$ denote the subalgebra of $\Xi$ that is generated
by $\{\c1,\c2,\c3,\d1,\d2\,\}$. 
Let $[\Xi_0]$ denote the subspace of $\Xi$ 
that is spanned by 
$\{\;\sym{\xi}ij\;\vert\; \xi \in \Xi_0;\;i,j\geq0\;\}$.
\end{definition}

\begin{remark}
\label{remnoc1}
Equations (\ref{sdrel}) and (\ref{symrel1}) imply that 
$\sym{z^kw}ij = 0$ if $w$ begins or ends with $\d1$.
(\ref{cdcrel}), (\ref{ddcrel}), (\ref{dddrel}) 
and $\d1\c1 = \c1\d1 = d_1^2 = 0$ allow to eliminate all $\d1$-s in the
middle of words. Thus $[\Xi_0]$ is spanned by elements of the form 
$\sym{z^k\,w}ij$ where $w$ is a word in the letters $\c1,\c2,\d2$. 
\end{remark}

\noindent
Let us do some auxiliary calculations in $[\Xi_0]$ based on Proposition \ref{propsym}.
\begin{equation}
\sym{z^k\c2}ij \eq{symrel1} \frac12\sym{z^k\c2}ij  + \frac12 \sym{z^k\s1\c2\s1}ij  \eq{screl} 
\frac12\sym{z^k(\c2+\c3)}ij \e \frac12\sym{z^k(z-\c1)}ij = 0
\label{zc0rel}
\end{equation}
$z^2\d2\c1(1+\s1) \eq{zzdcrel} z^3\d2(1+\s1)$ implies 
$\sym{z^{k}\d2}01 
= \sym{z^{k+1}\d2}00 
= 0$, so by Proposition \ref{propsym}
\begin{equation}
\label{zzd0rel}
\sym{z^k\d2}ij = 0 \mbox{\quad for }\, k \geq 2.
\end{equation}

\noindent
We can generalise this a little more:
\begin{equation}
\label{zzdw0rel}
\sym{z^k\d2w}ij = 0 \mbox{\quad if }\, k \geq 2 \mbox{ and } w \mbox{ is
a word in the letters } \c1 \mbox { and }\c2.
\end{equation}
This is shown by induction on the length of $w$. If $w$ is empty, then the statement 
is (\ref{zzd0rel}), otherwise if $w$ begins with $\c2$ then $\d2w = 0$,
so let us assume $w=\c1u$.
$$
2\sym{z^k\d2\c1u}ij \eq{cdsrel} \sym{z^k\d2\c1(1+\s1)u}ij + \sym{z^k\d2\d1u}ij \eq{zzdcrel}
\sym{z^{k+1}\d2u}ij + \sym{z^{k+1}\d2\s1u}ij + \sym{z^k\d2\d1u}ij
$$
By (\ref{screl}) and (\ref{symrel1}), the second term is equal to 
$\sym{z^{k+1}\d2u^\prime}ij$ where $u^\prime = \s1u\s1$ is obtained by replacing in $u$ every $\c2$ 
by $z-\c1-\c2$. If $u$ is empty or begins with $\c1$ then $\sym{z^k\d2\d1u}ij = 0$,
otherwise $u=\c2v$ and $\sym{z^k\d2\d1u}ij \eq{ddcrel} \sym{z^{k+1}\d2\c1v}ij$.
So the induction hypothesis applies to all three terms on the right side of the 
equation above and we are done.
\medskip

\noindent
We are interested in $[\Xi_0]$ because it contains the element $\sym{y}13$ we want to calculate.
\begin{eqnarray*}
\sym{y}{i+1}{j+1} - \sym{z\d2}{i+1}{j+1} & = & \sym{\c1y\c1 - z\c1\d2\c1}ij 
\eq{yrel} - \sym{z\c1\s1\s2\d1\c2}ij
\eq{screl} - \sym{z\s1\s2\c3\d1\c2}ij \\
&\eq{symrel1}& -\sym{z^2\s1\s2\d1\c2}ij + \sym{z\s2\c2\d1\c2}ij 
\eqz{zzssdcrel}{cdsrel} -\sym{z^2\d2\c1}ij + \sym{z(\c2-\d2)\d1\c2}ij  \\
&\eqz{cdcrel}{ddcrel} &  -\sym{z^2\d2}i{j+1} + \sym{z\c1\d2\c1}ij - \sym{z^2\d2\c1}ij  
\eq{zzd0rel} \sym{z\d2}{i+1}{j+1}.
\end{eqnarray*}
This allows us to express $\sym{y}ij$ as element of $[\Xi_0]$:
\begin{eqnarray}
\label{yzdrel}
\sym{y}ij &=& 2\sym{z\d2}ij \mbox{ for  }i,j \geq 1
\end{eqnarray}

\noindent
Next, we show that for $i,\,j,\,k \geq 0$ and arbitrary $\xi \in \Xi$, we have 
\begin{eqnarray}
\label{endred2}
2\sym{\xi\,c_1^2\c2}ij & = & \sym{z\xi}i{j+2} - \sym{\xi}i{j+3} \\
\label{endred3}
2\sym{\xi\,\d2\c1\c2}ij & = & 2\sym{z\xi\,\d2}i{j+1} - \sym{\xi\,\d2}i{j+2} \\
\label{endred4}
2\sym{\xi\,\c2\c1\c2}ij & = & \sym{z\xi\,\c2}i{j+1} - \sym{\xi\,\c2}i{j+2} + \sym{\xi\,\c1\d2}i{j+1}.
\end{eqnarray}
Equations (\ref{endred2})-(\ref{endred4}) are shown simultaneously:
$$
\sym{u(\c1\!-\!\d1)\c2}ij  \eq{cdsrel} \sym{u\c1\s1\c2}ij 
\eqz{screl}{symrel1} \sym{u\c1\c3}ij  \e \sym{u\c1(z\!-\!\c1\!-\!\c2)}ij 
\e \sym{zu}i{j+1} - \sym{u}i{j+2} -\sym{u\c1\c2}ij 
$$
This implies 
\quad $2\sym{u\,\c1\c2}ij  \e \sym{zu}i{j+1} - \sym{u}i{j+2} +\sym{u\,\d1\c2}ij.$\quad
If we set $u = \xi\c1$, $u = \xi\d2$, $u = \xi\c2$ and apply $\c1\d1 \!\e\! 0$,
$\d2\d1\c2 \!\eq{ddcrel}\! z\d2\c1$, $\c2\d1\c2 \!\eq{cdcrel}\! \c1\d2\c1$ we obtain
(\ref{endred2}), (\ref{endred3}), (\ref{endred4}), respectively.
\medskip

\begin{remark}
\label{remtyp}
Using the relations we have found so far, one can show
that $[\Xi_0]$ is spanned by the elements of 
the following forms:
\begin{itemize}
\item
$\sym{z^kc_2^2}{i}{i+1}$ with $k \geq 1$, $i \geq 0$,
\item
$\sym{z^k\d2}ij$ with $k \in \{0,1\}$ and $0 \leq i < j$,
\item
$\sym{\d2c_1^{2n}\d2}ij$ with $0 \leq i < j$ and $n \geq 1$. 
\end{itemize}
\end{remark}

\noindent
We only give a few hints how this can be shown.
\newcommand{\sig}[1]{\sigma_{#1}^{}}
Let $\sig{n} := \frac12(1+\s1)c_1^n$, then $\sig{n}\sig{m} = \sig{n+m}$
and $\c1\sig{n} = \sig{n}\c1 = \sig{n+1}$. 
By (\ref{cdsrel}) and Remark \ref{remnoc1}, $[\Xi_0]$ then is spanned
by elements of the form $\sym{z^k\sig{0}w_1^{}\sig{n_1}w_2^{}\cdots\sig{n_{l-1}}w_l^{}\sig{0}}pq$
with $n_i \geq 1$ and $w_i^{} = d_2$ or $w_i^{} = c_2^j$; this will be called
{\sl normal form} from now on. 
The $w_i$ are called {\sl segments} and we define the simplicity of
an element in normal form by: $k$ + number of $\d2$-segments.
For a linear combination $x$ of words of the same simplicity
we say $x \sim 0$ if any element $\sym{z^kuxw}pq$ in normal form 
can be expressed by a sum of elements
with simpler normal forms (i.e.~elements that contain more $z$-s or more $d_2$-s).
Of course, by $x \sim y$ we mean $x-y \sim 0$.
It is not hard to show that 
{
\newcounter{abccounter}
\renewcommand{\theequation}{\stepcounter{abccounter}\alph{abccounter}\addtocounter{equation}{-1}}
\begin{eqnarray}
\sig{n}\c2\sig{m} &\sim& -\frac12\sig{n+m+1} \\
\sig{n}c_2^{i+3}\sig{m} &\si{cccrel}& -2\sig{n}c_2^{i+1}\sig1\c2\sig{m} \;\;\sinoref{a}\;\;
\sig{n}c_2^{i+1}\sig{m+2} \\
\sig{n}c_2^2\sig{m+1} + \sig{n+1}c_2^2\sig{m}&\sim& \frac54\sig{n+m+3}\\
\sig{n}c_2^2\sig1c_2^2\sig{m} &\sim& -\frac12\sig{n}c_2^5\sig{m}\\
\sig{n}c_2^2\sig1\d2\sig{m} &\si{cccdrel}& -2\sig{n}\c2\sig2\d2\sig{m}  \;\;\sinoref{a}\;\; \sig{n+3}\d2\sig{m}. 
\end{eqnarray}
}
Using (a) and (b), we may eliminate all $c_2^j$-segments with $j \neq 2$.
Due to (c), we can reduce the $\sig{n}$-s between two $c_2^2$-segments or between a $c_2^2$-segment
and a $\d2$-segment to become $\sig1$ and then apply (d) or (e). Like this
we are left only with two types of normal forms: the ones that contain only $\d2$-segments
and those which consist of a single $c_2^2$-segment.
In the first case we may eliminate $\d2\sig{2n-1}\d2$ 
and $z\d2\sig{2n}\d2$ 
by (\ref{dcccdrel}) and (\ref{dcdrel}).
If a normal form contains $\d2\sig{2n}\d2\sig{2m}\d2$ then after rotating the 
even ladders with (LS), we may apply the relation (LL) to transform it
into $\frac14z^{n+m+2}\d2$. In the second case we may use (c) to put
the $c_2^2$ in the middle to get $\sym{z^kc_2^2}ii = 0$ or $\sym{z^kc_2^2}i{i+1}$.
If $k=0$ we may apply equation (\ref{ccrel}).
This shows that any element of $[\Xi_0]$ is a linear combination of
elements of the types in Remark \ref{remtyp}.

\begin{remark}
Using the relations that will be found in the next section, 
the statement of Remark \ref{remtyp} 
still can be improved a little bit:
$[\Xi_0]$ is spanned by 
$\sym{z^kc_2^2}{i-1}{i}$, $\sym{\d2}0i$, $\sym{z\d2}0i$,
$\sym{z\d2}1i$, $\sym{z\d2}{i}{i+1}$ with $i,k \geq 1$.
\end{remark}

\subsection{Additional Relations in $[\Xi_0]$}   
\label{addrelsect}
Let us consider a family $(\H{n})_{n\geq1}$ of elements of $\Xi$, where
$H_n$ is the sum of the four possible ways of attaching the lower ends of a
$n$-ladder to the first or second strand of $1$, and the upper ends to the
third strand. For example, $\H{3}$ looks like this:
\cepsfig{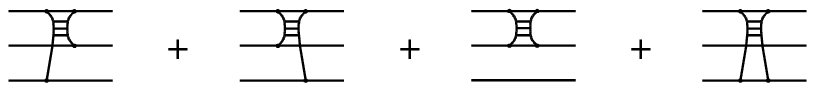}

\noindent
Let $R_n :=  \sym{H_n}01$.
Applying Lemma \ref{lemmab} a) twice, we see that $\H{n}\c1 = \c1\H{n}$.
Obviously $\mu\cdot\H{n} = \H{n}$, so 
$(1-\mu)\cdot\H{n}\c1 \e \H{n}\c1 - \c1\H{n} = 0$, thus 
\begin{eqnarray}
\label{R0rel}
R_n \e 0
\end{eqnarray}
We now present another way to write $R_n$ as element of $[\Xi_0]$
in order to obtain new relations. 
Let $\H{n,i}$ denote the $i$-th term of $\H{n}$ in the order of the
picture above, then  
$\H{n,2} = \s1\H{n,1}\s1$ and $\H{n,4} = \s1\H{n,3}\s1$. Thus by (\ref{symrel1}) we
have 
$$R_n = 2\sym{\H{n,1}}01 + 2\sym{\H{n,3}}01.$$
$\H{n,3}$ contains a $(n+2)$-ladder, so by (LS) and (LIHX) 
$$\H{n,3} = \left\{
\begin{array}{ll}
c_2^{n+2} & \mbox{if $n$ is even} \\
-c_2^{n+2} -\tau_3\cdot c_2^{n+2} & \mbox{if $n$ is odd}
\end{array}\right\}$$ 
and by (\ref{symrel1}) we get
$$2\sym{\H{n,3}}01 = \big(3\cdot(-1)^n-1\big)\,\sym{c_2^{n+2}}01.$$ 
Lemma \ref{lemmab} a) implies 
\cepsfig{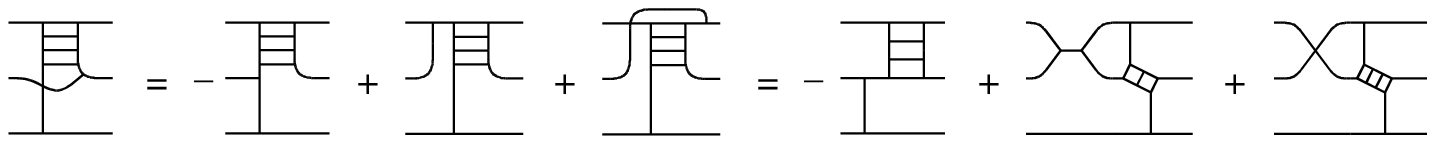}
\begin{equation}
\H{n,1} = -\c1\H{n-1,3}+\d2\,\l{n}+\s2\,\l{n+1} \mbox{\quad for }n\geq 1.
\label{H1rel}
\end{equation}

\noindent
Now $\sym{\c1\H{n-1,3}}01 = \sym{\H{n-1,3}}11 = 0$ because of $\mu\cdot\H{n-1,3} = \H{n-1,3}$.
The second and third term of $2\sym{\H{n,1}}01$ are according to (\ref{H1rel})
\begin{eqnarray*}
2\sym{\d2\l{n}}01 &\eqq{symrel1}& \sym{\d2\l{n}(1+\s1)}01 \eq{lsrel} 
\sym{\d2\d1\c2c_1^{n-1} + \d2\c3c_1^n}01 
\eq{ddcrel} \sym{z\d2c_1^n + \d2(z-\c1-\c2)c_1^n}01 \\
&=& 2\sym{z\d2c_1^n}01 - \sym{\d2c_1^{n+1}}01 \e 2\sym{z\d2}0{n+1}-\sym{\d2}0{n+2}  \\
2\sym{\s2\l{n+1}}01 & \eqq{slrel}&  2\sym{-\l{n+1}+c_2^{n}\c1\d2 + c_2^{n+1}\c3}01
\e 2\sym{c_2^{n}\c1\d2}01 + 2\sym{zc_2^{n+1}}01 -2\sym{c_2^{n+1}}02 -2\sym{c_2^{n+2}}01. 
\end{eqnarray*}

\noindent
The first equation is valid only for $n \geq 2$; for $n=1$ we use $\l1 = \c1\c2$.
In the last equation we used $\sym{\l{n}}ij = 0$, which is a consequence of
$\tau_3\cdot\l{n} = -\l{n}$ (relation (AS)). Summarising the calculations of this subsection
so far, we may state that
\begin{eqnarray}
\label{R1exp}
R_1 &=& 2\sym{\d2\c1\c2}01 +2\sym{\c2\c1\d2}01 
+2\sym{zc_2^2}01 -2\sym{c_2^2}02 -6\sym{c_2^3}01 \mbox{\quad and}\\
\label{Rnexp}
R_n &=& 2\sym{z\d2}0{n+1} - \sym{\d2}0{n+2} +2\sym{c_2^n\c1\d2}01 
+2\sym{zc_2^{n+1}}01 -2\sym{c_2^{n+1}}02 +3\big((-1)^n-1)\big)\sym{c_2^{n+2}}01.\mbox{\quad\quad}
\end{eqnarray}
with $n\geq 2$ are trivial elements of $[\Xi_0]$. In the spirit of Remark 
\ref{remtyp} one is able to present $R_n$ as linear combination of
simple elements. We need this for $n \leq 3$ which requires a 
lengthy computation that is done in the appendix. The results are
\begin{eqnarray}
\label{r1rel}
R_1 & = & 3\sym{\d2}12 - 6\sym{zc_2^2}01\\
\label{r2rel}
R_2 & = & 2\sym{\d2c_1^2\d2}01 - 3\sym{zd_2}03\\
\label{r3rel}
R_3 & = & 4\sym{\d2c_1^2\d2}02 - 6\sym{zd_2}04. 
\end{eqnarray}
To simplify the calculation of $R_2$ and $R_3$, the following formula has been used
\begin{equation}
\label{drel}
\sym{\d2}ij \e 2\sym{zc_2^2}{i-1}{j-1} \mbox{\quad for }\; i,j \geq 1.
\end{equation}
Because of $\sym{\c1\d2\c1 - 2zc_2^2}01 \eq{r1rel} \frac13R_1 \eq{R0rel} 0$ 
Proposition \ref{propsym} implies (\ref{drel}).

\noindent
Now we finally may combine (\ref{r2rel}), (\ref{r3rel}) and (\ref{yzdrel}) 
to obtain the desired result:
\begin{eqnarray}
\nonumber
0 & \eq{R0rel} & \frac13 R_3  - \frac23(\c1R_2+R_2\c1) \eq{symrel3} 
 \frac43\sym{\d2c_1^2\d2}02 - 2\sym{zd_2}04 
- \frac43\sym{\d2c_1^2\d2}11 + 2\sym{zd_2}13 \\
&& - \frac43\sym{\d2c_1^2\d2}02 + 2\sym{zd_2}04 \e
 2\sym{zd_2}13 \eq{yzdrel} \sym{y}13
\label{cycccrel}
\end{eqnarray}

\subsection{Calculations for parts 2 and 3}
Now, having the relation $\sym{y}13 \e 0$ in hand, we may 
stop messing around in $\Xi$ and proceed with 
the proof of the square-tunnelling
relation. But now that we have gone through so much trouble,
we first take profit out of our current knowledge to provide some more 
equations that are urgently needed in \cite{CC2} and \cite{CC3}. 

\noindent
We start with the observation that $\tau_3\cdot\c2\s2\c1\c2 \e \c2\s2\c2\c1 - \c2\s2\c1\c2$:
\cepsfig{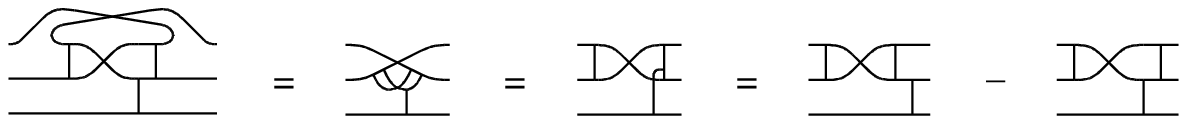}

\noindent
We obtain for any $i,j \geq 0$
\begin{eqnarray}
\nonumber
\sym{c_2^2}i{j+1} &=& \sym{\c2\s2\c2\c1}ij \e \sym{(1+\tau_3)\cdot\c2\s2\c1\c2}ij
\eqz{symrel1}{cdsrel} 2\sym{(\c2-\d2)\c1\c2}ij 
\eqz{endred4}{endred3} \sym{z\c2}i{j+1} - \sym{\c2}i{j+2} \\
\label{ccrel}
&&+\sym{\c1\d2}i{j+1}  -2\sym{z\d2}i{j+1} + \sym{\d2}i{j+2}     
\eq{zc0rel} \sym{\d2}{i+1}{j+1} + \sym{\d2}i{j+2}  - 2\sym{z\d2}i{j+1} 
\end{eqnarray}

\noindent
We use this for $i=1$ and $j=1$ to get an expression
for $c_1^2c_2^2c_1^2$:
\begin{eqnarray*}
0 &=& \c1 \big(\sym{c_2^2}12 - \sym{\d2}22 - \sym{\d2}13 + 2\sym{z\d2}12\big)
\eq{yzdrel} \c1 \big(\sym{c_2^2}12 - \sym{\d2}13 + \sym{y}12\big) \\
&=& (1+\tau_1)(1+\tau_3)(1+\tau_5)\cdot\big( 
c_1^2c_2^2c_1^2 - c_1^3c_2^2\c1 
- c_1^2\d2c_1^3 + c_1^4\d2\c1
+ c_1^2yc_1^2 - c_1^3y\c1 \big)
\end{eqnarray*}
Using $\tau_1\cdot c_1^2\xi \e c_1^2\xi$ and 
$\tau_3\cdot c_1^2c_2^2c_1^2 \e
\tau_5\cdot c_1^2c_2^2c_1^2 \e c_1^2c_2^2c_1^2$,  we obtain 
\begin{equation}
c_1^2c_2^2c_1^2 \e
\frac14(1+\tau_3)(1+\tau_5)\cdot\big( 
c_1^3y\c1 - c_1^4\d2\c1 
+ c_1^3c_2^2\c1 
- c_1^2yc_1^2 
+ c_1^2\d2c_1^3 
\big).
\end{equation}

\noindent
Another relation that will be needed in \cite{CC3} is the following: ($n\geq 1$)
\begin{eqnarray}
\nonumber
\d1\H{n,1} &\eq{H1rel}& -\d1\c1\H{n-1,3} + \d1\d2\l{n} + \d1\s2\l{n+1}\\
&\eqz{slrel}{screl}&  \d1\d2\l{n} - \d1\l{n+1} + \d1c_2^n\c1\d2 +\d1c_2^{n+1}\c1\s2
\end{eqnarray}

\noindent
Equations (\ref{symrel2}) and (\ref{cycccrel}) imply $\sym{y}n{n+2} = 0$ in $\ti F(6)$. In \cite{CC2} we need a similar relation in $F(6)$. 
This means that the whole calculation should be repeated in $F(6)$, 
i.e.~using the following 
relation instead of (x):
\cepsfig{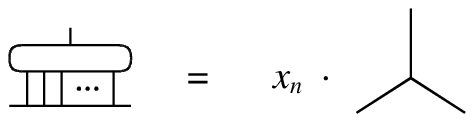}

\noindent
Interpreting $\sym{y}ij$ as element of $F(6)$ we would then obtain 
a relation of the following form (note that 
only $x_n$-s with $n \leq 4$ can occur and that $x_1 = 2t$, $x_2 = t^2$,
$x_4 = \frac43t\,x_3 - \frac13t^4$, see \cite{Vo}):
$$\sym{y}n{n+2} \e \sum_{i} \lambda_i\,\xi_i \mbox{\quad with }\lambda_i \in \Q[t,x_3]
- \Q
\mbox{ and } \xi_i \in F(6)$$
For $n=2$ we may even take the relation into a form in which 
all $\xi_i$ are of the form $\sym{y}**$ (which is not possible for $n=1$):
\begin{eqnarray}
\nonumber
3\sym{y}24 &=& 9t\,\sym{y}23 + 3t\,\sym{y}14 -9t^2\,\sym{y}13
+6t^2\,\sym{y}04 + (4t^3+2x_3)\,\sym{y}12 \\
&&-18t^3\,\sym{y}03
+4t(4t^3-x_3)\,\sym{y}02 +8t^2(x_3-t^3)\,\sym{y}01
\label{f6rel}
\end{eqnarray}

\subsection{The Square-Tunnelling Relation}

First, we use $\tau_1\cdot \c1yc_1^3 \e \s1\c1yc_1^3 \eq{yrel} \c1yc_1^3$ and 
$\tau_5\cdot\c1yc_1^3 \e \c1yc_1^3\s1 \eq{cdsrel} \c1yc_1^3$ to
reformulate equation (\ref{cycccrel}):
\begin{equation}
\label{therel}
(1+\tau_3)\cdot(\c1yc_1^3- c_1^3y\c1) \e
\frac14(1-\mu)(1+\tau_1)(1+\tau_3)(1+\tau_5)\cdot \c1yc_1^3 \e \frac14\sym{y}13  \eq{cycccrel} 0
\end{equation}

\newcommand{\Di}[2]{D_{#2,#1}}

\noindent
We turn our attention towards $\ti F(n+4)$. For $n \geq 0$ let 
$\Di1n^{}$, $\Di1n^\prime$, $\Di2n^{}$ and $\Di2n^\prime$ denote
the elements of $\ti F(n+4)$ that are represented by:
\cepsfig{}

\noindent
Furthermore let $D_n^{} := \Di1n^{}- \Di2n^{}$ and $D_n^{\prime} := \Di1n^\prime- \Di2n^\prime$. 

\begin{lemma}
\label{lemmswap}
$D_n^{} = 0$ for all $n \geq 0$ .
\end{lemma}
\begin{proof}
The assertion is true for $n=0$, since $\Di10^{} = \Di20^{}$. We proceed by induction and
assume that $n \geq 1$ and  $D_i^{} = 0$ for $i < n$.
We glue a line with $n$ legs with labels from $1$ to $n$ between the ends $3$ and $4$
of $\c1yc_1^3$ and rename the ends $1,2,4,6$ to $n+1,n+2,n+3,n+4$; the result is
$\Di1n^{}$: 
\cepsfig{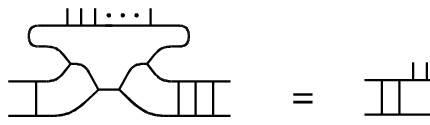}

\noindent
If we do the same with $\tau_3\cdot\c1yc_1^3$ and apply (AS)
at the $n$ legs, we get $(-1)^n\Di1n^{\prime}$.

\noindent
A similar statement can be done for $\c1yc_1^3$, $\tau_3\cdot\c1yc_1^3$
and $\Di2n^{}$, $(-1)^n\Di2n^{\prime}$, respectively, thus equation (\ref{therel}) implies:
\begin{equation}
\label{DDrel1}
D_n^{} + (-1)^nD_n^\prime \e \Di1n^{}  -\Di2n^{} + (-1)^n\Di1n^\prime - (-1)^n\Di2n^\prime \e 0. 
\end{equation}
By permuting neighbouring ends of both ladders, we get the following representation of $\Di1n^\prime$:
\cepsfig{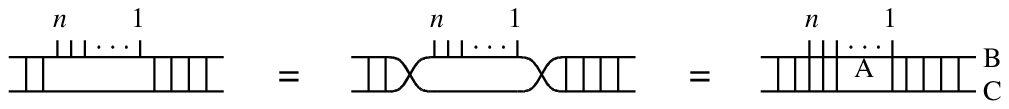}

\noindent
By Lemma \ref{lemmab} a) we can push all $n$ legs off the lower strand, one after another.
Let us push them all to the right, then we get $3^n$ terms which have legs 
at the edges indicated by the letters A, B, C in the picture above. The coefficient
of each term is $(-1)^{\mbox{\footnotesize number of legs in position A}}$.

We do the same trick for $\Di2n^\prime$ and get a similar sum with the 2-ladder
and the 4-ladder exchanged. If there are legs in position B or C, then there
must be less than $n$ legs between the two ladders. By induction assumption all these
terms in $\Di1n^\prime$ are equal to those in $\Di2n^\prime$. 
So in $D_n^\prime = \Di1n^\prime - \Di2n^\prime$ only the two terms having all $n$ legs in 
position A remain. The process reverses the order of the legs, so the 
remaining two diagrams are $\Di1n^{}$ and $\Di2n^{}$:
\begin{equation}
\label{DDrel2}
D_n^\prime \e  (-1)^n\Di1n^{} - (-1)^n\Di2n^{} \e (-1)^nD_n^{}
\end{equation}
The equations (\ref{DDrel1}) and (\ref{DDrel2}) imply the assertion of the lemma.
\end{proof}
\qed
\bigskip

\noindent
{\bf Proof of Theorem \ref{theo2}}\quad
We consider a graph containing two ladders with disjoint interiors.
Let us remove the interiors of these ladders and call the 
components of the remaining graph that contain ends of both ladders
{\sl joining}. A sequence of different consecutive
edges from one ladder to the other is called {\sl joining path}.
Remember that the condition
of Theorem \ref{theo2} is satisfied iff
\begin{itemize}
\item 
either there are two or more joining components, or
\item
there is a single joining component that contains two ends of one of the ladders or 
\item 
there is a single joining component that contains a circle.
\end{itemize}
The last case splits up into two subcases: either there are at least
two different paths joining the ladders or there is a unique joining path.
So to prove Theorem \ref{theo2} we have to show that in the following four situations we may 
exchange the square and the $4$-ladder:
$$
\begin{array}{ccccc}
\mbox{a)} & \hcenterfig{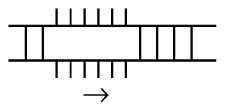} && 
\mbox{b)} & \hcenterfig{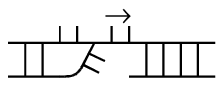} \\[0.3cm]
\mbox{c)} & \hcenterfig{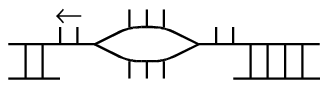} && 
\mbox{d)} & \hcenterfig{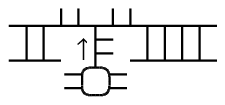}
\end{array}
$$
In situation a), we choose a joining path $p$ and push all edges arriving
at $p$ by Lemma \ref{lemmab} a) through one of the ladders, as indicated by
the little arrow. We can do the same for the graph in which both ladders 
have been exchanged. Both times we get a linear
combination of diagrams in which $p$ has become a single edge.
The only difference is that in each term the ladders have been 
exchanged. By Lemma \ref{lemmswap} the two expressions
are identical.

In situation b) we have two joining paths $p_1$ and $p_2$ that meet each other
and then coincide. We first push away 
all edges arriving at the common part of $p_1$ and $p_2$, to get
diagrams in which $p_1$ and $p_2$ have only an edge $e$ in common 
($e$ is the end of one of the ladders $L$). By Lemma \ref{lemmswap} a) we may express each of
these diagram by the three graphs we obtain if we cut off $p_2$ from $e$ and glue
it to one of the three other ends of $L$. So we may express any graph in situation
b) by a sum of graphs of situation a).

Exactly the same trick allows to reduce case c) to b). Finally in situation
d), we have a joining path $p_1$, a circle $l$ and a path $p_2$ connecting
$l$ and $p_1$. We assume that $p_2$ is a single edge (otherwise proceed
as above to clean $p_2$ as indicated by the arrow) and apply (IHX) at $p_2$. The result is a difference
of two graphs of type c).
This completes the proof of Theorem \ref{theo2}.
\qed
\bigskip

\noindent
{\bf Proof of Theorem \ref{coro}}\quad
The univalent vertices in the following pictures are marked
by bullets. If a graph in $B^u$ has weakly connected ladders of lengths $2$ and $4$,
it typically looks like this (the omitted parts of the graph are supposed to
be located inside the boxes):
\cepsfig{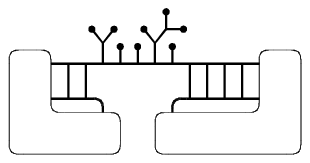}

\noindent
Since all univalent vertices carry the same colour, the (AS)-relation implies:
\eqnfig{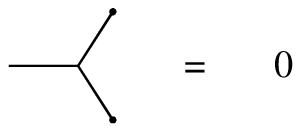}{\star}

\noindent
Lemma \ref{lemmab} a) yields
\cepsfig{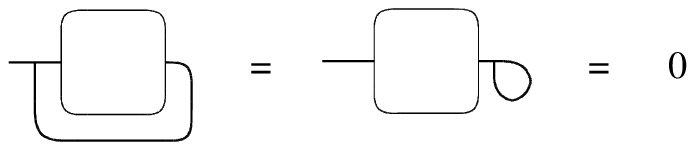}

\noindent
if the box contains no univalent vertices. So we assume that there 
is at least one univalent vertex inside each box and consider diagrams that
look like this:
\cepsfig{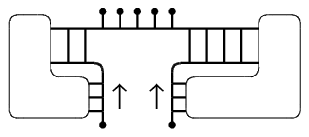}

\noindent
As before, we push all disturbing edges one by one through 
the ladders to obtain graphs that are trivial because of ($\star$), graphs that 
are of type b) or d) (in the proof of Theorem \ref{theo2}) and graphs that look like this:
\cepsfig{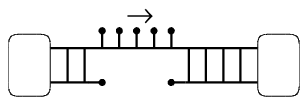}

\noindent
In the latter case, we push away all legs between the ladders.
Because of ($\star$) we only have to consider the terms in which all 
these legs end up in the two boxes. 
Thus to prove the theorem, it remains to show that we may exchange ladders in 
graphs of this type:
\cepsfig{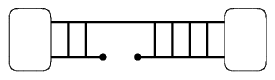}

\noindent
This is done by gluing univalent vertices of the same colour to the ends $3$ and $4$ in
equation (\ref{therel}) so that $\xi$ and $\tau_3\cdot\xi$ become the same element. 

\noindent
So in the case of unicoloured univalent vertices we may drop the
somehow artificial condition of strong connectivity and state
that the square-tunnelling relation holds in all connected graphs.
\qed

\begin{verbatim}
 
 
e-mail:jan@kneissler.info
http://www.kneissler.info
 
\end{verbatim}

%
%
%
\newpage
\setcounter{page}{1}
\renewcommand{\thepage}{A\arabic{page}}
{
\newcommand{\symb}[3]{[#1]_{#2}^{#3}}
\newcommand{\zeroind}[1]{{0_{\mathrm #1}}}
\newcommand{\neweq}[0]{\vspace{0.2cm} \\}

\section*{Appendix}
Here we give a detailed calculation that allows to express 
the elements $R_1$, $R_2$ and $R_3$  (given by (\ref{R1exp}) and (\ref{Rnexp}) in section \ref{addrelsect}) 
as linear combination of basic elements (according to Remark \ref{remtyp}). 
In every step, we use one of the following equations
($i,j \geq 0$ and $u,v$ arbitrary words):

\begin{displaymath}
\begin{array}{rrclcc}
\hspace{0.8cm}&\symb{u\,c_1^2c_2^{}}ij & = & \frac12\symb{zu}i{j+2}-\frac12\symb{u}i{j+3} &
\hspace{0.8cm} & \mbox{(A)} \neweq 
&\symb{u\,d_2^{}c_1^{}c_2^{}}ij & = & \symb{zu\,d_2^{}}i{j+1}-\frac12\symb{u\,d_2^{}}i{j+2}  && \mbox{(B)}  \neweq
&\symb{u\,c_2^{}c_1^{}c_2^{}}ij & = & \frac12\symb{zu\,c_2^{}}i{j+1}-\frac12\symb{u\,c_2^{}}i{j+2}+\frac12\symb{u\,c_1^{}d_2^{}}i{j+1} && \mbox{(C)} \neweq
&\symb{u\,c_2^3\,v}ij & = & \symb{zu\,c_2^2\,v}ij-2\symb{u\,c_2^{}c_1^{}c_2^{}\,v}ij+\symb{u\,d_2^{}c_1^{}c_2^{}\,v}ij+\symb{u\,c_2^{}c_1^{}d_2^{}\,v}ij && \mbox{(D)} \neweq
&\symb{u\,c_2^2c_1^{}d_2^{}\,v}ij & = & 2\symb{zu\,c_2^{}c_1^{}d_2^{}\,v}ij-2\symb{u\,c_2^{}c_1^2d_2^{}\,v}ij+\symb{u\,d_2^{}c_1^2d_2^{}\,v}ij && \mbox{(E)} \neweq
&\symb{u\,c_2^{}}i{j+2} & = & \symb{zu\,c_2^{}}i{j+1}-2\symb{u\,c_2^{}c_1^{}c_2^{}}ij+\symb{u\,c_1^{}d_2^{}}i{j+1} && \mbox{(C$^\prime$)} \neweq
&\symb{u\,c_2^{}c_1^{}c_2^{}\,v}ij & = & \frac12\symb{zu\,c_2^2\,v}ij-\frac12\symb{u\,c_2^3\,v}ij +\frac12\symb{u\,d_2^{}c_1^{}c_2^{}\,v}ij+\frac12\symb{u\,c_2^{}c_1^{}d_2^{}\,v}ij && \mbox{(D$^\prime$)} \neweq
&\symb{z\,c_2^2}ij & = & \frac12\symb{d_2^{}}{i+1}{j+1} && \mbox{(F)}
\end{array}
\end{displaymath}
The equations (A)-(F) are (\ref{endred2}), (\ref{endred3}), (\ref{endred4}), 
(\ref{cccrel}), (\ref{cccdrel}), (\ref{drel});
(C$^\prime$) and (D$^\prime$) are equivalent to (C) and (D). 
Note that equation (F) is not used in the calculation of $R_1$, because
the value of $R_1$ is an essential ingredient in the proof of (\ref{drel}).
\medskip

\noindent
Elements of the form $\symb{z^k\,w}ij$ (where $w$ is a word in the letters 
$c_1^{},\,c_2^{},\,d_2^{}$) are trivial in each of the following cases:
(see (\ref{dcdrel}), (\ref{zc0rel}), (\ref{zzdw0rel}), (\ref{dcccdrel}))

{
\begin{enumerate}
\setlength{\itemsep}{0.05cm}
\item[a)] 
$w$ contains $d_2^{}c_2^{}$ or $c_2^{}d_2^{}$ or $d_2^{}c_1^{}d_2^{}$ as subword,
\item[b)] 
$w = c_2^{}$,
\item[c)] 
$k \geq 2$ and $w=d_2^{}$ or $w=d_2^{}c_1c_2^{n}$,
\item[d)]
$w = d_2^{}c_1^{}c_2^{}c_1^{}d_2^{}$,
\item[e)] 
$k \geq 1$ and $w = d_2^{}c_1^{2}d_2^{}$,
\item[f)]
$w$ is a palindrome
and $i = j$.
\end{enumerate}
}

\noindent
The basic elements are named $e_1,\, \ldots,\, e_k$.
In the calculations each term $t$ of the left side of equations (A) to (F) 
has one of the following three properties:
{\begin{itemize}
\setlength{\itemsep}{0.0cm}
\item
$t$ or $\mu\cdot t$ is one of the basic elements, or
\item
$t$ or $\mu\cdot t$ belongs to one of the classes a) to f) and is henceforth 0, or
\item
$t$ or $\mu\cdot t$ has been calculated in a previous equation.
\end{itemize}}
\noindent
We indicate the second case by writing $\zeroind{x}$ (x one of the letters a to f),
and the third case by the equation number in brackets. 
If $\mu\cdot t$ (the reversed word) instead of $t$ satisfies the condition, the
sign of the coefficient of $t$ has to be changed according to (\ref{symrel1}):
$$\symb{w}ij \;\;=\;\; -\mu\cdot\symb{w}ij \;\;=\;\; -\symb{\mu\cdot w}ji.$$ 

\noindent
In each step the result is a linear combination of basic elements
which is written as $k$-dimensional row-vector.
}
%
%
%
{
\newpage
%
%
%

\newcommand{\symb}[3]{[#1]_{#2}^{#3}}
\newcommand{\equstack}[1]{{\stackrel{\mbox{\tiny#1}}{=}}}
\newcommand{\zeroind}[1]{{0_{\mathrm #1}}}
\newcommand{\vecmin}[0]{\mbox{-}}
\newcommand{\nextline}[0]{\vspace{0.1cm} \\}
\newcommand{\neweq}[0]{\vspace{0.2cm} \\}
\newcommand{\zsep}[0]{\,}

\subsection*{Calculation of $R_1$}

\noindent
Basic elements: 
$\;\; e_1\,=\,\symb{z\zsep c_2^{2}}01$,
$\;\; e_2\,=\,\symb{d_2^{}}03$,
$\;\; e_3\,=\,\symb{d_2^{}}12$,
$\;\; e_4\,=\,\symb{z\zsep d_2^{}}02$.

\begin{displaymath}
\begin{array}{rclccr}
 \symb{d_2^{}c_1^{2}c_2^{}}00 & \equstack{(A)} & 
\multicolumn{3}{l}{
\frac{1}{2}\symb{z\zsep d_2^{}}02-\frac{1}{2}\symb{d_2^{}}03} & \nextline
 & = & \frac{1}{2}e_{4}-\frac{1}{2}e_{2} & = &
(0,\, \vecmin \frac{1}{2},\, 0,\, \frac{1}{2}) & (1)\neweq
 \symb{z\zsep d_2^{}c_1^{}c_2^{}}00 & \equstack{(B)} & 
\multicolumn{3}{l}{
\symb{z^2\zsep d_2^{}}01-\frac{1}{2}\symb{z\zsep d_2^{}}02} & \nextline
 & = & \zeroind{c}-\frac{1}{2}e_{4} & = &
(0,\, 0,\, 0,\, \vecmin \frac{1}{2}) & (2)\neweq
 \symb{c_2^{2}c_1^{}d_2^{}}00 & \equstack{(E)} & 
\multicolumn{3}{l}{
2\symb{z\zsep c_2^{}c_1^{}d_2^{}}00-2\symb{c_2^{}c_1^{2}d_2^{}}00+\symb{d_2^{}c_1^{2}d_2^{}}00} & \nextline
 & = & -2(2)+2(1)+\zeroind{f} & = &
(0,\, \vecmin 1,\, 0,\, 2) & (3)\neweq
 \symb{c_2^{}c_1^{}c_2^{}}01 & \equstack{(C)} & 
\multicolumn{3}{l}{
\frac{1}{2}\symb{z\zsep c_2^{}}02-\frac{1}{2}\symb{c_2^{}}03+\frac{1}{2}\symb{d_2^{}}12} & \nextline
 & = & \zeroind{b}+\zeroind{b}+\frac{1}{2}e_{3} & = &
(0,\, 0,\, \frac{1}{2},\, 0) & (4)\neweq
 \symb{c_2^{2}c_1^{}c_2^{}}00 & \equstack{(D$^\prime$)} & 
\multicolumn{3}{l}{
\frac{1}{2}\symb{z\zsep c_2^{3}}00-\frac{1}{2}\symb{c_2^{4}}00+\frac{1}{2}\symb{c_2^{}d_2^{}c_1^{}c_2^{}}00+\frac{1}{2}\symb{c_2^{2}c_1^{}d_2^{}}00} & \nextline
 & = & \zeroind{f}+\zeroind{f}+\zeroind{a}+\frac{1}{2}(3) & = &
(0,\, \vecmin \frac{1}{2},\, 0,\, 1) & (5)\neweq
 \symb{d_2^{}c_1^{}c_2^{}}01 & \equstack{(B)} & 
\multicolumn{3}{l}{
\symb{z\zsep d_2^{}}02-\frac{1}{2}\symb{d_2^{}}03} & \nextline
 & = & e_{4}-\frac{1}{2}e_{2} & = &
(0,\, \vecmin \frac{1}{2},\, 0,\, 1) & (6)\neweq
 \symb{d_2^{}c_1^{}c_2^{}}10 & \equstack{(B)} & 
\multicolumn{3}{l}{
\symb{z\zsep d_2^{}}11-\frac{1}{2}\symb{d_2^{}}12} & \nextline
 & = & \zeroind{f}-\frac{1}{2}e_{3} & = &
(0,\, 0,\, \vecmin \frac{1}{2},\, 0) & (7)\neweq
 \symb{c_2^{2}}02 & \equstack{(C$^\prime$)} & 
\multicolumn{3}{l}{
\symb{z\zsep c_2^{2}}01-2\symb{c_2^{2}c_1^{}c_2^{}}00+\symb{c_2^{}c_1^{}d_2^{}}01} & \nextline
 & = & e_{1}-2(5)-(7) & = &
(1,\, 1,\, \frac{1}{2},\, \vecmin 2) & (8)\neweq
 \symb{c_2^{3}}01 & \equstack{(D)} & 
\multicolumn{3}{l}{
\symb{z\zsep c_2^{2}}01-2\symb{c_2^{}c_1^{}c_2^{}}01+\symb{d_2^{}c_1^{}c_2^{}}01+\symb{c_2^{}c_1^{}d_2^{}}01} & \nextline
 & = & e_{1}-2(4)+(6)-(7) & = &
(1,\, \vecmin \frac{1}{2},\, \vecmin \frac{1}{2},\, 1) & (9)\neweq
\end{array}
\end{displaymath}

\begin{displaymath}
\begin{array}{rcl}
R_1 & = & 2\symb{d_2^{}c_1^{}c_2^{}}01+2\symb{c_2^{}c_1^{}d_2^{}}01
+2\symb{z\,c_2^2}01 -2\symb{c_2^2}02 -6\symb{c_2^3}01 \nextline
& = & 2(6)-2(7)+2e_1-2(8)-6(9)
\;\;=\;\; (\vecmin 6,\, 0,\, 3,\, 0) 
\end{array}
\end{displaymath}

\newpage
\subsection*{Calculation of $R_2$}

\noindent
We use equation (F) at one place (line (9)). The basic elements are: 
$$e_1\,=\,\symb{z^2\zsep c_2^{2}}01,
 \;\; e_2\,=\,\symb{d_2^{}}04,
 \;\; e_3\,=\,\symb{z\zsep d_2^{}}03,
 \;\; e_4\,=\,\symb{z\zsep d_2^{}}12,
 \;\; e_5\,=\,\symb{d_2^{}c_1^{2}d_2^{}}01.$$

\begin{displaymath}
\begin{array}{rclccr}
 \symb{z\zsep d_2^{}c_1^{2}c_2^{}}00 & \equstack{(A)} & 
\multicolumn{3}{l}{
\frac{1}{2}\symb{z^2\zsep d_2^{}}02-\frac{1}{2}\symb{z\zsep d_2^{}}03} & \nextline
 & = & \zeroind{c}-\frac{1}{2}e_{3} & = &
(0,\, 0,\, \vecmin \frac{1}{2},\, 0,\, 0) & (1)\neweq
 \symb{z\zsep c_2^{2}c_1^{}d_2^{}}00 & \equstack{(E)} & 
\multicolumn{3}{l}{
2\symb{z^2\zsep c_2^{}c_1^{}d_2^{}}00-2\symb{z\zsep c_2^{}c_1^{2}d_2^{}}00+\symb{z\zsep d_2^{}c_1^{2}d_2^{}}00} & \nextline
 & = & \zeroind{c}+2(1)+\zeroind{e} & = &
(0,\, 0,\, \vecmin 1,\, 0,\, 0) & (2)\neweq
 \symb{z\zsep c_2^{2}c_1^{}c_2^{}}00 & \equstack{(D$^\prime$)} & 
\multicolumn{3}{l}{
\frac{1}{2}\symb{z^2\zsep c_2^{3}}00-\frac{1}{2}\symb{z\zsep c_2^{4}}00+\frac{1}{2}\symb{z\zsep c_2^{}d_2^{}c_1^{}c_2^{}}00+\frac{1}{2}\symb{z\zsep c_2^{2}c_1^{}d_2^{}}00} & \nextline
 & = & \zeroind{f}+\zeroind{f}+\zeroind{a}+\frac{1}{2}(2) & = &
(0,\, 0,\, \vecmin \frac{1}{2},\, 0,\, 0) & (3)\neweq
 \symb{d_2^{}c_1^{}c_2^{}}20 & \equstack{(B)} & 
\multicolumn{3}{l}{
\symb{z\zsep d_2^{}}21-\frac{1}{2}\symb{d_2^{}}22} & \nextline
 & = & -e_{4}+\zeroind{f} & = &
(0,\, 0,\, 0,\, \vecmin 1,\, 0) & (4)\neweq
 \symb{d_2^{}c_1^{}c_2^{}}02 & \equstack{(B)} & 
\multicolumn{3}{l}{
\symb{z\zsep d_2^{}}03-\frac{1}{2}\symb{d_2^{}}04} & \nextline
 & = & e_{3}-\frac{1}{2}e_{2} & = &
(0,\, \vecmin \frac{1}{2},\, 1,\, 0,\, 0) & (5)\neweq
 \symb{z\zsep d_2^{}c_1^{}c_2^{}}01 & \equstack{(B)} & 
\multicolumn{3}{l}{
\symb{z^2\zsep d_2^{}}02-\frac{1}{2}\symb{z\zsep d_2^{}}03} & \nextline
 & = & \zeroind{c}-\frac{1}{2}e_{3} & = &
(0,\, 0,\, \vecmin \frac{1}{2},\, 0,\, 0) & (6)\neweq
 \symb{z\zsep c_2^{}c_1^{}c_2^{}}01 & \equstack{(C)} & 
\multicolumn{3}{l}{
\frac{1}{2}\symb{z^2\zsep c_2^{}}02-\frac{1}{2}\symb{z\zsep c_2^{}}03+\frac{1}{2}\symb{z\zsep d_2^{}}12} & \nextline
 & = & \zeroind{b}+\zeroind{b}+\frac{1}{2}e_{4} & = &
(0,\, 0,\, 0,\, \frac{1}{2},\, 0) & (7)\neweq
 \symb{z\zsep d_2^{}c_1^{}c_2^{}}10 & \equstack{(B)} & 
\multicolumn{3}{l}{
\symb{z^2\zsep d_2^{}}11-\frac{1}{2}\symb{z\zsep d_2^{}}12} & \nextline
 & = & \zeroind{c}-\frac{1}{2}e_{4} & = &
(0,\, 0,\, 0,\, \vecmin \frac{1}{2},\, 0) & (8)\neweq
\frac12\symb{d_2^{}}13 \;\;\equstack{(F)}\;\;
 \symb{z\zsep c_2^{2}}02 & \equstack{(C$^\prime$)} & 
\multicolumn{3}{l}{
\symb{z^2\zsep c_2^{2}}01-2\symb{z\zsep c_2^{2}c_1^{}c_2^{}}00+\symb{z\zsep c_2^{}c_1^{}d_2^{}}01} & \nextline
 & = & e_{1}-2(3)-(8) & = &
(1,\, 0,\, 1,\, \frac{1}{2},\, 0) & (9)\neweq
 \symb{c_2^{}c_1^{}c_2^{}}02 & \equstack{(C)} & 
\multicolumn{3}{l}{
\frac{1}{2}\symb{z\zsep c_2^{}}03-\frac{1}{2}\symb{c_2^{}}04+\frac{1}{2}\symb{d_2^{}}13} & \nextline
 & = & \zeroind{b}+\zeroind{b}+(9) & = &
(1,\, 0,\, 1,\, \frac{1}{2},\, 0) & (10)\neweq
 \symb{d_2^{}c_1^{2}c_2^{}}10 & \equstack{(A)} & 
\multicolumn{3}{l}{
\frac{1}{2}\symb{z\zsep d_2^{}}12-\frac{1}{2}\symb{d_2^{}}13} & \nextline
 & = & \frac{1}{2}e_{4}-(9) & = &
(\vecmin 1,\, 0,\, \vecmin 1,\, 0,\, 0) & (11)\neweq
 \symb{c_2^{2}c_1^{}d_2^{}}01 & \equstack{(E)} & 
\multicolumn{3}{l}{
2\symb{z\zsep c_2^{}c_1^{}d_2^{}}01-2\symb{c_2^{}c_1^{2}d_2^{}}01+\symb{d_2^{}c_1^{2}d_2^{}}01} & \nextline
 & = & -2(8)+2(11)+e_{5} & = &
(\vecmin 2,\, 0,\, \vecmin 2,\, 1,\, 1) & (12)\neweq
 \symb{z\zsep c_2^{3}}01 & \equstack{(D)} & 
\multicolumn{3}{l}{
\symb{z^2\zsep c_2^{2}}01-2\symb{z\zsep c_2^{}c_1^{}c_2^{}}01+\symb{z\zsep d_2^{}c_1^{}c_2^{}}01+\symb{z\zsep c_2^{}c_1^{}d_2^{}}01} & \nextline
 & = & e_{1}-2(7)+(6)-(8) & = &
(1,\, 0,\, \vecmin \frac{1}{2},\, \vecmin \frac{1}{2},\, 0) & (13)\neweq
 \symb{c_2^{3}}02 & \equstack{(D)} & 
\multicolumn{3}{l}{
\symb{z\zsep c_2^{2}}02-2\symb{c_2^{}c_1^{}c_2^{}}02+\symb{d_2^{}c_1^{}c_2^{}}02+\symb{c_2^{}c_1^{}d_2^{}}02} & \nextline
 & = & (9)-2(10)+(5)-(4) & = &
(\vecmin 1,\, \vecmin \frac{1}{2},\, 0,\, \frac{1}{2},\, 0) & (14)\neweq
\end{array}
\end{displaymath}

\begin{displaymath}
\begin{array}{rcl}
R_2 & = & 2\symb{z\,d_2^{}}03-\symb{d_2^{}}04
+2\symb{c_2^2c_1^{}d_2^{}}01
+2\symb{z\,c_2^3}01 -2\symb{c_2^3}02 \nextline
& = & 2e_3-e_2+2(12)+2(13)-2(14)
\;\;=\;\; (0,\, 0,\, \vecmin3,\, 0,\, 2) 
\end{array}
\end{displaymath}

\newpage
\subsection*{Calculation of $R_3$}

\noindent
Relation (F) will be used in line (24); furthermore it allows to 
use the basic element $e_1$ in two different forms:
$$    e_1\,=\,\symb{z\zsep c_2^{2}}12\,=\,\frac12\symb{d_2^{}}23,
 \;\;\, e_2\,=\,\symb{z^3\zsep c_2^{2}}01,
 \;\;\, e_3\,=\,\symb{d_2^{}}05,$$
$$ e_4\,=\,\symb{z\zsep d_2^{}}04,
 \;\;\;\;\;\, e_5\,=\,\;\symb{z\zsep d_2^{}}13,
 \;\;\;\;\;\, e_6\,=\,\symb{d_2^{}c_1^{2}d_2^{}}02.$$

\begin{displaymath}
\begin{array}{rclccr}
 \symb{d_2^{}c_1^{}c_2^{}}03 & \equstack{(B)} & 
\multicolumn{3}{l}{
\symb{z\zsep d_2^{}}04-\frac{1}{2}\symb{d_2^{}}05} & \nextline
 & = & e_{4}-\frac{1}{2}e_{3} & = &
(0,\, 0,\, \vecmin \frac{1}{2},\, 1,\, 0,\, 0) & (1)\neweq
 \symb{d_2^{}c_1^{}c_2^{}}21 & \equstack{(B)} & 
\multicolumn{3}{l}{
\symb{z\zsep d_2^{}}22-\frac{1}{2}\symb{d_2^{}}23} & \nextline
 & = & \zeroind{f}-e_{1} & = &
(\vecmin 1,\, 0,\, 0,\, 0,\, 0,\, 0) & (2)\neweq
 \symb{c_2^{}c_1^{}c_2^{}}12 & \equstack{(C)} & 
\multicolumn{3}{l}{
\frac{1}{2}\symb{z\zsep c_2^{}}13-\frac{1}{2}\symb{c_2^{}}14+\frac{1}{2}\symb{d_2^{}}23} & \nextline
 & = & \zeroind{b}+\zeroind{b}+e_{1} & = &
(1,\, 0,\, 0,\, 0,\, 0,\, 0) & (3)\neweq
 \symb{z^2\zsep c_2^{2}c_1^{}c_2^{}}00 & \equstack{(D$^\prime$)} & 
\multicolumn{3}{l}{
\frac{1}{2}\symb{z^3\zsep c_2^{3}}00-\frac{1}{2}\symb{z^2\zsep c_2^{4}}00+\frac{1}{2}\symb{z^2\zsep c_2^{}d_2^{}c_1^{}c_2^{}}00+\frac{1}{2}\symb{z^2\zsep c_2^{2}c_1^{}d_2^{}}00} & \nextline
 & = & \zeroind{f}+\zeroind{f}+\zeroind{a}+\zeroind{c} & = &
(0,\, 0,\, 0,\, 0,\, 0,\, 0) & (4)\neweq
 \symb{z\zsep d_2^{}c_1^{2}c_2^{}}10 & \equstack{(A)} & 
\multicolumn{3}{l}{
\frac{1}{2}\symb{z^2\zsep d_2^{}}12-\frac{1}{2}\symb{z\zsep d_2^{}}13} & \nextline
 & = & \zeroind{c}-\frac{1}{2}e_{5} & = &
(0,\, 0,\, 0,\, 0,\, \vecmin \frac{1}{2},\, 0) & (5)\neweq
 \symb{z^2\zsep c_2^{}c_1^{}c_2^{}}01 & \equstack{(C)} & 
\multicolumn{3}{l}{
\frac{1}{2}\symb{z^3\zsep c_2^{}}02-\frac{1}{2}\symb{z^2\zsep c_2^{}}03+\frac{1}{2}\symb{z^2\zsep d_2^{}}12} & \nextline
 & = & \zeroind{b}+\zeroind{b}+\zeroind{c} & = &
(0,\, 0,\, 0,\, 0,\, 0,\, 0) & (6)\neweq
 \symb{z\zsep d_2^{}c_1^{}c_2^{}}11 & \equstack{(B)} & 
\multicolumn{3}{l}{
\symb{z^2\zsep d_2^{}}12-\frac{1}{2}\symb{z\zsep d_2^{}}13} & \nextline
 & = & \zeroind{c}-\frac{1}{2}e_{5} & = &
(0,\, 0,\, 0,\, 0,\, \vecmin \frac{1}{2},\, 0) & (7)\neweq
 \symb{z\zsep d_2^{}c_1^{2}c_2^{}}01 & \equstack{(A)} & 
\multicolumn{3}{l}{
\frac{1}{2}\symb{z^2\zsep d_2^{}}03-\frac{1}{2}\symb{z\zsep d_2^{}}04} & \nextline
 & = & \zeroind{c}-\frac{1}{2}e_{4} & = &
(0,\, 0,\, 0,\, \vecmin \frac{1}{2},\, 0,\, 0) & (8)\neweq
 \symb{z\zsep d_2^{}c_1^{}c_2^{}}20 & \equstack{(B)} & 
\multicolumn{3}{l}{
\symb{z^2\zsep d_2^{}}21-\frac{1}{2}\symb{z\zsep d_2^{}}22} & \nextline
 & = & \zeroind{c}+\zeroind{f} & = &
(0,\, 0,\, 0,\, 0,\, 0,\, 0) & (9)\neweq
 \symb{z\zsep c_2^{}c_1^{}c_2^{}}02 & \equstack{(C)} & 
\multicolumn{3}{l}{
\frac{1}{2}\symb{z^2\zsep c_2^{}}03-\frac{1}{2}\symb{z\zsep c_2^{}}04+\frac{1}{2}\symb{z\zsep d_2^{}}13} & \nextline
 & = & \zeroind{b}+\zeroind{b}+\frac{1}{2}e_{5} & = &
(0,\, 0,\, 0,\, 0,\, \frac{1}{2},\, 0) & (10)\neweq
 \symb{z^2\zsep c_2^{2}}02 & \equstack{(C$^\prime$)} & 
\multicolumn{3}{l}{
\symb{z^3\zsep c_2^{2}}01-2\symb{z^2\zsep c_2^{2}c_1^{}c_2^{}}00+\symb{z^2\zsep c_2^{}c_1^{}d_2^{}}01} & \nextline
 & = & e_{2}-2(4)+\zeroind{c} & = &
(0,\, 1,\, 0,\, 0,\, 0,\, 0) & (11)\neweq
 \symb{d_2^{}c_1^{2}c_2^{}}02 & \equstack{(A)} & 
\multicolumn{3}{l}{
\frac{1}{2}\symb{z\zsep d_2^{}}04-\frac{1}{2}\symb{d_2^{}}05} & \nextline
 & = & \frac{1}{2}e_{4}-\frac{1}{2}e_{3} & = &
(0,\, 0,\, \vecmin \frac{1}{2},\, \frac{1}{2},\, 0,\, 0) & (12)\neweq
\end{array}
\end{displaymath}
\newpage
\begin{displaymath}
\begin{array}{rclccr}
 \symb{z\zsep d_2^{}c_1^{}c_2^{}}02 & \equstack{(B)} & 
\multicolumn{3}{l}{
\symb{z^2\zsep d_2^{}}03-\frac{1}{2}\symb{z\zsep d_2^{}}04} & \nextline
 & = & \zeroind{c}-\frac{1}{2}e_{4} & = &
(0,\, 0,\, 0,\, \vecmin \frac{1}{2},\, 0,\, 0) & (13)\neweq
 \symb{d_2^{}c_1^{}c_2^{}c_1^{}c_2^{}}01 & \equstack{(C)} & 
\multicolumn{3}{l}{
\frac{1}{2}\symb{z\zsep d_2^{}c_1^{}c_2^{}}02-\frac{1}{2}\symb{d_2^{}c_1^{}c_2^{}}03+\frac{1}{2}\symb{d_2^{}c_1^{2}d_2^{}}02} & \nextline
 & = & \frac{1}{2}(13)-\frac{1}{2}(1)+\frac{1}{2}e_{6} & = &
(0,\, 0,\, \frac{1}{4},\, \vecmin \frac{3}{4},\, 0,\, \frac{1}{2}) & (14)\neweq
 \symb{c_2^{2}c_1^{}d_2^{}}20 & \equstack{(E)} & 
\multicolumn{3}{l}{
2\symb{z\zsep c_2^{}c_1^{}d_2^{}}20-2\symb{c_2^{}c_1^{2}d_2^{}}20+\symb{d_2^{}c_1^{2}d_2^{}}20} & \nextline
 & = & -2(13)+2(12)-e_{6} & = &
(0,\, 0,\, \vecmin 1,\, 2,\, 0,\, \vecmin 1) & (15)\neweq
 \symb{z\zsep c_2^{3}}02 & \equstack{(D)} & 
\multicolumn{3}{l}{
\symb{z^2\zsep c_2^{2}}02-2\symb{z\zsep c_2^{}c_1^{}c_2^{}}02+\symb{z\zsep d_2^{}c_1^{}c_2^{}}02+\symb{z\zsep c_2^{}c_1^{}d_2^{}}02} & \nextline
 & = & (11)-2(10)+(13)-(9) & = &
(0,\, 1,\, 0,\, \vecmin \frac{1}{2},\, \vecmin 1,\, 0) & (16)\neweq
 \symb{z\zsep c_2^{2}c_1^{}d_2^{}}10 & \equstack{(E)} & 
\multicolumn{3}{l}{
2\symb{z^2\zsep c_2^{}c_1^{}d_2^{}}10-2\symb{z\zsep c_2^{}c_1^{2}d_2^{}}10+\symb{z\zsep d_2^{}c_1^{2}d_2^{}}10} & \nextline
 & = & \zeroind{c}+2(8)+\zeroind{e} & = &
(0,\, 0,\, 0,\, \vecmin 1,\, 0,\, 0) & (17)\neweq
 \symb{d_2^{}c_1^{}c_2^{3}}01 & \equstack{(D)} & 
\multicolumn{3}{l}{
\symb{z\zsep d_2^{}c_1^{}c_2^{2}}01-2\symb{d_2^{}c_1^{}c_2^{}c_1^{}c_2^{}}01+\symb{d_2^{}c_1^{}d_2^{}c_1^{}c_2^{}}01+\symb{d_2^{}c_1^{}c_2^{}c_1^{}d_2^{}}01} & \nextline
 & = & -(17)-2(14)+\zeroind{a}+\zeroind{d} & = &
(0,\, 0,\, \vecmin \frac{1}{2},\, \frac{5}{2},\, 0,\, \vecmin 1) & (18)\neweq
 \symb{z\zsep c_2^{2}c_1^{}c_2^{}}10 & \equstack{(C)} & 
\multicolumn{3}{l}{
\frac{1}{2}\symb{z^2\zsep c_2^{2}}11-\frac{1}{2}\symb{z\zsep c_2^{2}}12+\frac{1}{2}\symb{z\zsep c_2^{}c_1^{}d_2^{}}11} & \nextline
 & = & \zeroind{f}-\frac{1}{2}e_{1}-\frac{1}{2}(7) & = &
(\vecmin \frac{1}{2},\, 0,\, 0,\, 0,\, \frac{1}{4},\, 0) & (19)\neweq
 \symb{z^2\zsep c_2^{3}}01 & \equstack{(D)} & 
\multicolumn{3}{l}{
\symb{z^3\zsep c_2^{2}}01-2\symb{z^2\zsep c_2^{}c_1^{}c_2^{}}01+\symb{z^2\zsep d_2^{}c_1^{}c_2^{}}01+\symb{z^2\zsep c_2^{}c_1^{}d_2^{}}01} & \nextline
 & = & e_{2}-2(6)+\zeroind{c}+\zeroind{c} & = &
(0,\, 1,\, 0,\, 0,\, 0,\, 0) & (20)\neweq
 \symb{z\zsep c_2^{2}c_1^{}d_2^{}}01 & \equstack{(E)} & 
\multicolumn{3}{l}{
2\symb{z^2\zsep c_2^{}c_1^{}d_2^{}}01-2\symb{z\zsep c_2^{}c_1^{2}d_2^{}}01+\symb{z\zsep d_2^{}c_1^{2}d_2^{}}01} & \nextline
 & = & \zeroind{c}+2(5)+\zeroind{e} & = &
(0,\, 0,\, 0,\, 0,\, \vecmin 1,\, 0) & (21)\neweq
 \symb{z\zsep c_2^{4}}01 & \equstack{(D)} & 
\multicolumn{3}{l}{
\symb{z^2\zsep c_2^{3}}01-2\symb{z\zsep c_2^{}c_1^{}c_2^{2}}01+\symb{z\zsep d_2^{}c_1^{}c_2^{2}}01+\symb{z\zsep c_2^{}c_1^{}d_2^{}c_2^{}}01} & \nextline
 & = & (20)+2(19)-(17)+\zeroind{a} & = &
(\vecmin 1,\, 1,\, 0,\, 1,\, \frac{1}{2},\, 0) & (22)\neweq
 \symb{z\zsep c_2^{2}c_1^{}c_2^{}}01 & \equstack{(D$^\prime$)} & 
\multicolumn{3}{l}{
\frac{1}{2}\symb{z^2\zsep c_2^{3}}01-\frac{1}{2}\symb{z\zsep c_2^{4}}01+\frac{1}{2}\symb{z\zsep c_2^{}d_2^{}c_1^{}c_2^{}}01+\frac{1}{2}\symb{z\zsep c_2^{2}c_1^{}d_2^{}}01} & \nextline
 & = & \frac{1}{2}(20)-\frac{1}{2}(22)+\zeroind{a}+\frac{1}{2}(21) & = &
(\frac{1}{2},\, 0,\, 0,\, \vecmin \frac{1}{2},\, \vecmin \frac{3}{4},\, 0) & (23)\neweq
%
\frac12\symb{d_2^{}}14 & \equstack{(F)} &
\multicolumn{3}{l}{
 \symb{z\zsep c_2^{2}}03 \;\;\; \equstack{(C$^\prime$)} \;\;\; 
\symb{z^2\zsep c_2^{2}}02-2\symb{z\zsep c_2^{2}c_1^{}c_2^{}}01+\symb{z\zsep c_2^{}c_1^{}d_2^{}}02} & \nextline
 & = & (11)-2(23)-(9) & = &
(\vecmin 1,\, 1,\, 0,\, 1,\, \frac{3}{2},\, 0) & (24)\neweq
 \symb{d_2^{}c_1^{2}c_2^{}}11 & \equstack{(A)} & 
\multicolumn{3}{l}{
\frac{1}{2}\symb{z\zsep d_2^{}}13-\frac{1}{2}\symb{d_2^{}}14} & \nextline
 & = & \frac{1}{2}e_{5}-(24) & = &
(1,\, \vecmin 1,\, 0,\, \vecmin 1,\, \vecmin 1,\, 0) & (25)\neweq
 \symb{d_2^{}c_1^{}c_2^{}}12 & \equstack{(B)} & 
\multicolumn{3}{l}{
\symb{z\zsep d_2^{}}13-\frac{1}{2}\symb{d_2^{}}14} & \nextline
 & = & e_{5}-(24) & = &
(1,\, \vecmin 1,\, 0,\, \vecmin 1,\, \vecmin \frac{1}{2},\, 0) & (26)\neweq
 \symb{c_2^{2}c_1^{}d_2^{}}11 & \equstack{(E)} & 
\multicolumn{3}{l}{
2\symb{z\zsep c_2^{}c_1^{}d_2^{}}11-2\symb{c_2^{}c_1^{2}d_2^{}}11+\symb{d_2^{}c_1^{2}d_2^{}}11} & \nextline
 & = & -2(7)+2(25)+\zeroind{f} & = &
(2,\, \vecmin 2,\, 0,\, \vecmin 2,\, \vecmin 1,\, 0) & (27)\neweq
\end{array}
\end{displaymath}
\newpage
\begin{displaymath}
\begin{array}{rclccr}
 \symb{c_2^{3}}12 & \equstack{(D)} & 
\multicolumn{3}{l}{
\symb{z\zsep c_2^{2}}12-2\symb{c_2^{}c_1^{}c_2^{}}12+\symb{d_2^{}c_1^{}c_2^{}}12+\symb{c_2^{}c_1^{}d_2^{}}12} & \nextline
 & = & e_{1}-2(3)+(26)-(2) & = &
(1,\, \vecmin 1,\, 0,\, \vecmin 1,\, \vecmin \frac{1}{2},\, 0) & (28)\neweq
 \symb{c_2^{3}c_1^{}c_2^{}}10 & \equstack{(C)} & 
\multicolumn{3}{l}{
\frac{1}{2}\symb{z\zsep c_2^{3}}11-\frac{1}{2}\symb{c_2^{3}}12+\frac{1}{2}\symb{c_2^{2}c_1^{}d_2^{}}11} & \nextline
 & = & \zeroind{f}-\frac{1}{2}(28)+\frac{1}{2}(27) & = &
(\frac{1}{2},\, \vecmin \frac{1}{2},\, 0,\, \vecmin \frac{1}{2},\, \vecmin \frac{1}{4},\, 0) & (29)\neweq
 \symb{c_2^{2}c_1^{}c_2^{}}20 & \equstack{(C)} & 
\multicolumn{3}{l}{
\frac{1}{2}\symb{z\zsep c_2^{2}}21-\frac{1}{2}\symb{c_2^{2}}22+\frac{1}{2}\symb{c_2^{}c_1^{}d_2^{}}21} & \nextline
 & = & -\frac{1}{2}e_{1}+\zeroind{f}-\frac{1}{2}(26) & = &
(\vecmin 1,\, \frac{1}{2},\, 0,\, \frac{1}{2},\, \frac{1}{4},\, 0) & (30)\neweq
 \symb{d_2^{}c_1^{}c_2^{}c_1^{}c_2^{}}10 & \equstack{(C)} & 
\multicolumn{3}{l}{
\frac{1}{2}\symb{z\zsep d_2^{}c_1^{}c_2^{}}11-\frac{1}{2}\symb{d_2^{}c_1^{}c_2^{}}12+\frac{1}{2}\symb{d_2^{}c_1^{2}d_2^{}}11} & \nextline
 & = & \frac{1}{2}(7)-\frac{1}{2}(26)+\zeroind{f} & = &
(\vecmin \frac{1}{2},\, \frac{1}{2},\, 0,\, \frac{1}{2},\, 0,\, 0) & (31)\neweq
 \symb{c_2^{3}c_1^{}d_2^{}}01 & \equstack{(D)} & 
\multicolumn{3}{l}{
\symb{z\zsep c_2^{2}c_1^{}d_2^{}}01-2\symb{c_2^{}c_1^{}c_2^{}c_1^{}d_2^{}}01+\symb{d_2^{}c_1^{}c_2^{}c_1^{}d_2^{}}01+\symb{c_2^{}c_1^{}d_2^{}c_1^{}d_2^{}}01} & \nextline
 & = & (21)+2(31)+\zeroind{d}+\zeroind{a} & = &
(\vecmin 1,\, 1,\, 0,\, 1,\, \vecmin 1,\, 0) & (32)\neweq
 \symb{c_2^{4}}02 & \equstack{(D)} & 
\multicolumn{3}{l}{
\symb{z\zsep c_2^{3}}02-2\symb{c_2^{}c_1^{}c_2^{2}}02+\symb{d_2^{}c_1^{}c_2^{2}}02+\symb{c_2^{}c_1^{}d_2^{}c_2^{}}02} & \nextline
 & = & (16)+2(30)-(15)+\zeroind{a} & = &
(\vecmin 2,\, 2,\, 1,\, \vecmin \frac{3}{2},\, \vecmin \frac{1}{2},\, 1) & (33)\neweq
 \symb{c_2^{5}}01 & \equstack{(D)} & 
\multicolumn{3}{l}{
\symb{z\zsep c_2^{4}}01-2\symb{c_2^{}c_1^{}c_2^{3}}01+\symb{d_2^{}c_1^{}c_2^{3}}01+\symb{c_2^{}c_1^{}d_2^{}c_2^{2}}01} & \nextline
 & = & (22)+2(29)+(18)+\zeroind{a} & = &
(0,\, 0,\, \vecmin \frac{1}{2},\, \frac{5}{2},\, 0,\, \vecmin 1) & (34)\neweq
\end{array}
\end{displaymath}

\begin{displaymath}
\begin{array}{rcl}
R_3 & = & 2\symb{z\,d_2^{}}04-\symb{d_2^{}}05
+2\symb{c_2^3c_1^{}d_2^{}}01
+2\symb{z\,c_2^4}01 -2\symb{c_2^4}02 -6\symb{c_2^5}01 \nextline
& = & 2e_4 -e_3 +2(32)+2(22)-2(33)-6(34)
\;\;=\;\; (0,\, 0,\, 0,\, \vecmin 6,\, 0,\, 4) 
\end{array}
\end{displaymath}

}


\begin{thebibliography}{00}

\bibitem{BN} D.~Bar-Natan, {\it On the Vassiliev knot invariants,} Topology {\bf 34} (1995), 423-472.

\bibitem{CC2} J.~A.~Kneissler, {\it On Spaces of Connected Graphs II: Relations in $\Lambda$}, Jour.~of Knot Theory and its Ramif.~Vol.~{\bf 10}, No. 5 (2001), 667-674.

\bibitem{CC3} J.~A.~Kneissler, {\it On Spaces of Connected Graphs III: The ladder filtration}, Jour.~of Knot Theory and its Ramif.~Vol.~{\bf 10}, No. 5 (2001), 675-686.

\bibitem{Ko} M.~Kontsevich, {\it Vassiliev's knot invariants,} Adv.~in Sov.~
Math.,~{\bf 16(2)} (1993), 137-150.

\bibitem{Va} V.~A. Vassiliev, {\it Cohomology of knot spaces,} Theory of Singularities and its Applications
(ed. V.~I. Arnold), Advances in Soviet Math., {\bf 1} (1990) 23-69.

\bibitem{Vo} Pierre Vogel, {\it Algebraic structures on modules of diagrams,} Universit\'e Paris VII preprint, July 1995 (revised 1997).

\end{thebibliography}
\end{document}